\documentclass{amsart}
\usepackage{amssymb}
\usepackage{a4}
\begin{document}
%\date{Version 3v as of 7 March 2006 by PBG}
\newtheorem{theorem}{Theorem}[section]
\newtheorem{definition}[theorem]{Definition}
\newtheorem{lemma}[theorem]{Lemma}
\newtheorem{example}[theorem]{Example}
\newtheorem{remark}[theorem]{Remark}
\newtheorem{corollary}[theorem]{Corollary}
\def\Rank{\operatorname{Rank}}
\def\ffrac#1#2{{\textstyle\frac{#1}{#2}}}
\def\qedbox{\hbox{$\rlap{$\sqcap$}\sqcup$}}
\def\id{\operatorname{id}}\def\Tr{\operatorname{Tr}}
\makeatletter
\renewcommand{\theequation}{%
\thesection.\alph{equation}} \@addtoreset{equation}{section}
\makeatother
\def\mapright#1{\smash{\mathop{\longrightarrow}\limits\sp{#1}}}
\title{Algebraic theory of affine curvature tensors}
\author{N.  Bla{\v z}i{\'c} ${}^1$,  P. Gilkey, S. Nik\v cevi\'c, and U. Simon}
%\begin{address}{NB: Faculty of Mathematics, University of Beograd, Studentski Trg. 16, P.P. 550,
%11000 Beograd, Srbija i Crna Gora. Email: {\it blazicn@matf.bg.ac.yu}}\end{address}
\begin{address}{PG: Mathematics Department, University of Oregon,
Eugene Or 97403 USA.\newline Email: {\it gilkey@darkwing.uoregon.edu}}
\end{address}
\begin{address}{SN: Mathematical Institute, Sanu,
Knez Mihailova 35, p.p. 367,
11001 Belgrade,
Serbia and Montenegro.
\newline Email: {\it stanan@mi.sanu.ac.yu}}\end{address}
\begin{address}
{US: Institut f\"ur Mathematik, Technische Universit\"at
Berlin\\ Strasse des 17. Juni 135, D-10623 Berlin, Deutschland\\
Email: {\it simon@math.tu-berlin.de}}\end{address}
\begin{abstract} 
We use curvature decompositions to construct generating sets
for the space of algebraic curvature tensors and for the space of tensors with
the same symmetries as those of a torsion free, Ricci symmetric connection;
the latter naturally appear in relative hypersurface theory.
\end{abstract}
\keywords{Algebraic curvature tensors, affine curvature tensors,
\newline 2000 {\it Mathematics Subject Classification.} 53B20\newline
\rm ${}^1$ N. Bla{\v z}i{\'c} passed away Monday 10 October 2005. This article is dedicated to his memory.}
\maketitle

\section{Introduction}

Let $V$ be a real vector space of dimension $m$; to simplify the discussion, we shall
assume that $m\ge4$ henceforth; similar results hold in dimensions $m=2$ and $m=3$. In Section \ref{sect-2},
we discuss the space of curvature operators $\mathfrak{R}(V)\subset\otimes^2V^*\otimes\operatorname{End}(V)$.
These are operators with the same symmetries as those of the curvature operator of a torsion free connection on the tangent bundle
of a smooth manifold. One has that
$\mathcal{R}\in\mathfrak{R}(V)$ if and only if for all $x,y,z\in V$,
\begin{eqnarray}
&&\mathcal{R}(x,y)z=-\mathcal{R}(y,x)z\quad\text{and}\label{eqn-1.a}\\
&&\mathcal{R}(x,y)z+\mathcal{R}(y,z)x+\mathcal{R}(z,x)y=0\,.\label{eqn-1.b}
\end{eqnarray}
Equation (\ref{eqn-1.b}) is called the {\it first Bianchi identity}.
We have, see for example Strichartz \cite{S88}, that
$$\dim\mathcal{R}(V)=\textstyle\frac13m^2(m^2-1)\,.$$

In Section
\ref{sect-3}, we discuss the space of algebraic curvature tensors
$\mathfrak{a}(V)\subset\otimes^4V^*$. This is the space of tensors with the same symmetries as that of the
curvature tensor defined by the Levi-Civita connection of a pseudo-Riemannian metric; $A\in\mathfrak{a}(V)$ if
and only if for all
$x,y,z,w\in V$,
\begin{eqnarray}
&&A(x,y,z,w)=-A(y,x,z,w),\label{eqn-1.c}\\
&&A(x,y,z,w)+A(y,z,x,w)+A(z,x,y,w)=0,\label{eqn-1.d}\\
&&A(x,y,z,w)=A(z,w,x,y),\label{eqn-1.e}\\
&&A(x,y,z,w)=-A(x,y,w,z)\,.\label{eqn-1.f}
\end{eqnarray}
We shall show in Theorem \ref{thm-3.2} that  identities (\ref{eqn-1.e}) and (\ref{eqn-1.f}) are equivalent in the presence of
identities (\ref{eqn-1.c}) and (\ref{eqn-1.d}). One has, see for example Strichartz \cite{S88}, that:
$$\dim\{\mathfrak{a}(V)\}=\textstyle\frac{m^2(m^2-1)}{12}\,.$$

If $\mathcal{R}\in\mathfrak{R}(V)$, it is natural to consider the traces:
\begin{equation}\label{eqn-1.g}
\begin{array}{l}
\rho_{14}(\mathcal{R})(x,y):=\operatorname{Tr}\{z\rightarrow\mathcal{R}(z,x)y\},\\
\rho_{24}(\mathcal{R})(x,y):=\operatorname{Tr}\{z\rightarrow\mathcal{R}(x,z)y\},\vphantom{\vrule height 11pt}\\
\rho_{34}(\mathcal{R})(x,y):=\operatorname{Tr}\{z\rightarrow\mathcal{R}(x,y)z\}\,.\vphantom{\vrule height 11pt}
\end{array}\end{equation}
The identities of Equations (\ref{eqn-1.a}) and (\ref{eqn-1.b}) show that:
\begin{equation}\label{eqn-1.h}
\begin{array}{l}
\rho_{24}(\mathcal{R})=-\rho_{14}(\mathcal{R})\quad\text{and}\\
\rho_{34}(\mathcal{R})(x,y)=-\rho_{14}(\mathcal{R})(x,y)+\rho_{14}(\mathcal{R})(y,x)\,.\vphantom{\vrule height 11pt}
\end{array}\end{equation}
In Section
\ref{sect-4}, we discuss the affine curvature operators $\mathfrak{F}(V)\subset\mathfrak{R}(V)$.
These are the operators with the same symmetries as those of an affine connection without torsion;
$\mathcal{F}\in\mathfrak{F}(V)$ if and only if for all $x,y,z\in V$,
\begin{eqnarray}
&&\mathcal{F}(x,y)z=-\mathcal{F}(y,x)z,\label{eqn-1.i}\\
&&\mathcal{F}(x,y)z+\mathcal{F}(y,z)x+\mathcal{F}(z,x)y=0,\label{eqn-1.j}\\
 &&\rho_{14}(\mathcal{F})(x,y)=\rho_{14}(\mathcal{F})(y,x),\label{eqn-1.k}\\
&&\rho_{34}(\mathcal{F})=0\,.\label{eqn-1.L}
\end{eqnarray} 
By Equation (\ref{eqn-1.h}), Equations (\ref{eqn-1.k}) and (\ref{eqn-1.L}) are equivalent in the
presence of Equations (\ref{eqn-1.i}) and (\ref{eqn-1.j}); thus these are the symmetries of the
curvature operator of a torsion free, Ricci symmetric connection on the tangent bundle of a smooth
manifold. Such curvature operators appear naturally as curvature operators of the induced and of the
conormal connections in relative hypersurface theory; see \cite{SSV91}.

\par The natural structure group of the spaces 
$\mathfrak{R}(V)$, $\mathfrak{a}(V)$, and $\mathfrak{F}(V)$ is the general linear group
$GL(V)$. Let $O(V,\langle\cdot,\cdot\rangle)$ be the orthogonal group associated to a non-degenerate symmetric
bilinear form $\langle\cdot,\cdot\rangle\in S^2(V^*)$ of signature $(p,q)$ on $V$. We can use $\langle\cdot,\cdot\rangle$ to raise
and lower indices and define an $O(V,\langle\cdot,\cdot\rangle)$ equivariant identification between
$\otimes^4V^*$ and $\otimes^2V^*\otimes\operatorname{End}(V)$ by means of the identity:
\begin{equation}\label{eqn-1.m}
R(x,y,z,w)=\langle\mathcal{R}(x,y)z,w\rangle\,.
\end{equation}
We let 
$$
  \mathfrak{r}(V)\subset\otimes^4V^*,\quad
  \mathfrak{A}(V,\langle\cdot,\cdot\rangle)\subset\otimes^2V^*\otimes\operatorname{End}(V),\quad
  \mathfrak{f}(V,\langle\cdot,\cdot\rangle)\subset\otimes^4V^*
$$
 be the subspaces associated to $\mathfrak{R}(V)$,
$\mathfrak{a}(V)$, and $\mathfrak{F}(V)$, respectively; $R\in\mathfrak{r}(V)$ if and only if for all $x,y,z,w\in V$, one has
\begin{eqnarray*}
&&R(x,y,z,w)=-R(y,x,z,w),\\
&&R(x,y,z,w)+R(y,z,x,w)+R(z,x,y,w)=0\,.
\end{eqnarray*}
We have $\mathcal{A}\in\mathfrak{A}(V,\langle\cdot,\cdot\rangle)$ if  and only if for all $x,y,z,w\in V$ one has:
\begin{eqnarray*}
&&\mathcal{A}(x,y)=-\mathcal{A}(y,x),\\
&&\mathcal{A}(x,y)z+\mathcal{A}(y,z)x+\mathcal{A}(z,x)y=0,\\
&&\langle\mathcal{A}(x,y)z,w\rangle=\langle\mathcal{A}(z,w)x,y\rangle,\\
&&\langle\mathcal{A}(x,y)z,w\rangle=-\langle\mathcal{A}(x,y)w,z\rangle,
\end{eqnarray*}
the last two identities being equivalent in the presence of the first two. Finally
$F\in\mathfrak{f}(V,\langle\cdot,\cdot\rangle)$ if and only if for all $x,y,z,w\in V$ one has:
\begin{eqnarray}
&&F(x,y,z,w)=-F(y,z,x,w),\label{eqn-1.n}\\
&&F(x,y,z,w)+F(y,z,x,w)+F(z,x,y,w)=0,\label{eqn-1.o}\\
&&\rho_{14}(F)(x,y)=\rho_{14}(F)(y,x),\label{eqn-1.p}\\
&&(\operatorname{id}\otimes\operatorname{Tr})F=0\,.\label{eqn-1.q}
\end{eqnarray}
Again, identities (\ref{eqn-1.p}) and (\ref{eqn-1.q}) are equivalent given the identities of Equations
(\ref{eqn-1.n}) and (\ref{eqn-1.o}).

\par The spaces $\mathfrak{A}(V,\langle\cdot,\cdot\rangle)$ and
$\mathfrak{f}(V,\langle\cdot,\cdot\rangle)$ depend upon the choice of the inner product; the space $\mathfrak{r}(V)$ does not.
Thus it is convenient to keep the distinction between subspaces of $\otimes^2V^*\otimes\operatorname{End}(V)$ and $\otimes^4V^*$;
this will play a crucial role in the proof of Theorem \ref{thm-4.2}. The spaces $\mathfrak{R}(V)$, $\mathfrak{A}(V)$, and
$\mathfrak{F}(V)$ are subspaces of $\otimes^2V^*\otimes\operatorname{End}(V)$; elements of these
spaces will be denoted by $\mathcal{R}$, $\mathcal{A}$, and
$\mathcal{F}$, respectively, and are operator valued bilinear forms. The spaces
$\mathfrak{r}(V)$, $\mathfrak{a}(V)$, and $\mathfrak{f}(V)$ are subspaces of
$\otimes^4V^*$; elements of these spaces will be denoted by $R$, $A$, and $F$, respectively,  and are quadralinear forms. We
have the inclusions:
$$\begin{array}{lllll}
\mathfrak{A}(V,\langle\cdot,\cdot\rangle)&\subset&\mathfrak{F}(V)&\subset&\mathfrak{R}(V),\\
\mathfrak{a}(V)&\subset&\mathfrak{f}(V,\langle\cdot,\cdot\rangle)&\subset&\mathfrak{r}(V)\,.
\end{array}$$

Let $\{e_i\}$ be a basis for $V$. If $\psi\in\otimes^2V^*$ and if
$\Theta\in\otimes^4V^*$, set
$$\psi_{ij}:=\psi(e_i,e_j)\quad\text{and}\quad\Theta_{ijkl}:=\Theta(e_i,e_j,e_k,e_l)\,.$$
Let $\{e^i\}$ be the associated dual basis for $V^*$. Then
$$\psi=\textstyle\sum_{ij}\psi_{ij}e^i\otimes e^j\quad\text{and}\quad
  \Theta=\textstyle\sum_{ijkl}\psi_{ijkl}e^i\otimes e^j\otimes e^k\otimes e^l\,.$$
If $\langle\cdot,\cdot\rangle$ is a non-degenerate inner product on $V$, let 
\begin{equation}\label{eqn-1.r}
\Xi_{ij}:=\langle e_i,e_j\rangle\quad\text{and}\quad\textstyle\sum_j\Xi^{ij}\Xi_{jk}=\delta^i_k
\end{equation}
where $\delta$ is the Kronecker symbol. One then has:
$$
\textstyle\sum_{ij}\Xi^{ij}\langle x,e_i\rangle e_j=x\quad\text{and}\quad
  \Tr\{\psi\}=\textstyle\sum_{ij}\Xi^{ij}\psi_{ij}\,.
$$

 We shall decompose the natural
action of $GL(V)$ and of $O(V,\langle\cdot,\cdot\rangle)$ on the spaces $\mathfrak{R}(V)$, $\mathfrak{a}(V)$, and
$\mathfrak{F}(V)$ as the direct sum of irreducible modules and use these decompositions to exhibit generating sets
for these spaces and to derive other natural geometric properties.

\par Our motivation in this paper is to study affine curvature
operators; as already stated above, these are the curvature operators
which naturally appear as curvature operators of the induced and
of the conormal connections in relative hypersurface theory.
Moreover, in this situation, there naturally appears
a metric, the so called relative metric,  which permits us to
raise and lower indices. Our aim is a characterization of the
affine curvature operators, arising from torsion free and Ricci symmetric
connections, in the space of all curvature operators arising from torsion free
connections.
Via the decomposition results of Section
\ref{sect-4}, these are characterized by the vanishing of the component
$W_3$. We will study the geometric meaning of the various components
in this decomposition, at least in the case of relative hypersurfaces,
in a subsequent paper.

\section{Curvature Operators}\label{sect-2}

In this section, we study operators with the same symmetries as those of a torsion free connection on the tangent bundle of a
smooth  manifold.

\subsection{Geometric representability of curvature operators}\label{sect-2.1}
If $\nabla$ is a
torsion free connection on the tangent bundle of a smooth manifold $M$, let
$\mathcal{R}^\nabla$ be the associated curvature operator; if $P\in M$ and if $x,y,z\in T_PM$, then
$$\mathcal{R}_P^\nabla(x,y)z:=\left\{\nabla_x\nabla_y-\nabla_y\nabla_x-\nabla_{[x,y]}\right\}z\,.
$$
One then has $\mathcal{R}_P^\nabla\in\mathfrak{R}(T_PM)$ since the symmetries of Equations (\ref{eqn-1.a}) and (\ref{eqn-1.b})
hold. Conversely, every curvature operator is geometrically representable by an torsion free connection:
\begin{theorem}\label{thm-2.1}
Let $\mathcal{R}\in\mathfrak{R}(V)$ be given. Regard $V$ as a smooth manifold in its own right. Let $0$ be the origin
of $V$ and identify $T_0V=V$. Then there exists a torsion free connection $\nabla$ on $V$ so that
$\mathcal{R}_0^\nabla=\mathcal{R}$.
\end{theorem}

\begin{proof}Let $\mathcal{R}\in\mathfrak{R}(V)$. Expand $\mathcal{R}(e_i,e_j)e_k=\sum_lR_{ijk}{}^le_l$
relative to some basis $\{e_i\}$ for $V$. Let $\{x_i\}$ be the associated dual coordinates; if $e\in V$, then $e=\sum_ix_i(e)e_i$.
Define a connection $\nabla$ on $TV$ by setting
\begin{eqnarray*}
&&\nabla_{\partial_{x_a}}\partial_{x_b}:=\textstyle\sum_d\Gamma_{ab}{}^d\partial_{x_d}\quad\text{for}\quad
\Gamma_{ab}{}^d:=-\textstyle\frac13\sum_cx_c\{\mathcal{R}_{acb}{}^d+\mathcal{R}_{bca}{}^d\}\,.
\end{eqnarray*}
Since $\nabla_{\partial_{x_a}}\partial_{x_b}=\nabla_{\partial_{x_b}}\partial_{x_a}$, $\nabla$ is torsion
free. As $\Gamma(0)=0$,
\begin{eqnarray*}
&&\mathcal{R}_0^\nabla(\partial_{x_i},\partial_{x_j})\partial_{x_k}
=\textstyle\sum_l(\partial_{x_i}\Gamma_{jk}{}^l-\partial_{x_j}\Gamma_{ik}{}^l)\partial_{x_l}\\
&=&-\textstyle\frac13\sum_l\{\mathcal{R}_{jik}{}^l+\mathcal{R}_{kij}{}^\ell
   -\mathcal{R}_{ijk}{}^\ell-\mathcal{R}_{kji}{}^l\}\partial_{x_l}\\
&=&-\textstyle\frac13\sum_l\{-2\mathcal{R}_{ijk}{}^l+\mathcal{R}_{kij}{}^l
 +\mathcal{R}_{jki}{}^l\}\partial_{x_l}
 =\mathcal{R}_{ijk}{}^l\partial_{x_l}\,.
\end{eqnarray*}
This completes the proof of the desired result.\end{proof}

\subsection{The Jacobi operator}\label{sect-2.2} This operator is defined by setting:
$$\mathcal{J}_{\mathcal{R}}(x)y:=\mathcal{R}(y,x)x\,.$$
It plays a central role in the study of geodesic sprays. 
The following result is known in the context of Riemannian geometry; it extends easily to the more general setting.
\begin{lemma}\label{lem-2.2} Let $\mathcal{R}\in\mathfrak{R}(V)$. If $\mathcal{J}_{\mathcal{R}}=0$, then $\mathcal{R}=0$.
\end{lemma}

\begin{proof} $\mathcal{J}_{\mathcal{R}}(x)$ is quadratic in $x$. The associated bilinear form is given by
$$\mathcal{J}_{\mathcal{R}}(x,y):z\rightarrow\textstyle\frac12\{\partial_{\varepsilon}
  \mathcal{J}_{\mathcal{R}}(x+\varepsilon
y)\}z|_{\varepsilon=0}=\textstyle\frac12\{\mathcal{R}(z,x)y+\mathcal{R}(z,y)x\}\,.$$
If $\mathcal{J}_{\mathcal{R}}(x)=0$ for all $x\in V$, one has the additional curvature symmetry
$$\mathcal{R}(z,x)y+\mathcal{R}(z,y)x=0$$
for all $x,y,z\in V$. We compute:
\begin{eqnarray*}
0&=&\mathcal{R}(x,y)z+\mathcal{R}(y,z)x+\mathcal{R}(z,x)y\\
&=&\mathcal{R}(x,y)z-\mathcal{R}(y,x)z-\mathcal{R}(x,z)y\\
&=&\mathcal{R}(x,y)z+\mathcal{R}(x,y)z+\mathcal{R}(x,y)z\,.
\end{eqnarray*}
The Lemma now follows.\end{proof}

\subsection{The action of the general linear group on $\mathfrak{R}(V)$}\label{sect-2.3} This action is not irreducible, but
decomposes as the direct sum of irreducible modules. Let
\begin{equation}\label{eqn-2.a}
\mathfrak{U}(V):=\ker\{\rho_{14}\}\cap\mathfrak{R}(V)\,.
\end{equation}
The decomposition $V^*\otimes V^*=\Lambda^2(V^*)\oplus S^2(V^*)$ is a $GL(V)$ equivariant decomposition of $V^*\otimes V^*$
into irreducible modules; we let $\pi_a$ and $\pi_s$ be the appropriate projections where
\begin{equation}\label{eqn-2.b}
\pi_a(\psi)_{ij}:=\textstyle\frac12(\psi_{ij}-\psi_{ji})\quad\text{and}\quad
  \pi_s(\psi)_{ij}:=\textstyle\frac12(\psi_{ij}+\psi_{ji})\,.
\end{equation}
We may therefore decompose $\rho_{14}=\pi_a\circ\rho_{14}\oplus\pi_s\circ\rho_{14}$ where
$\rho_{14}$ is as defined in Equation (\ref{eqn-1.g}).
One has the following result of Strichartz \cite{S88}:

\begin{theorem}\label{thm-2.3}
The map $\rho_{14}$ defines a $GL(V)$ equivariant short exact sequence
$$0\rightarrow\mathfrak{U}(V)\rightarrow\mathfrak{R}(V)\mapright{\rho_{14}}
\Lambda^2(V^*)\oplus S^2(V^*)\rightarrow0$$
which is equivariantly split by the map $\sigma_{\pi_a\circ\rho_{14}}\oplus\sigma_{\pi_s\circ\rho_{14}}$ where
\begin{eqnarray*}
&&\sigma_{\pi_a\circ\rho_{14}}(\omega)(x,y)z=\textstyle\frac{-1}{1+m}\{2\omega(x,y)z+\omega(x,z)y-\omega(y,z)x\}
\text{ for }\omega\in\Lambda^2(V^*),\\
&&\sigma_{\pi_s\circ\rho_{14}}(\psi)(x,y)z=\textstyle\frac1{1-m}\{\psi(x,z)y-\psi(y,z)x\}
\text{ for }\psi\in S^2(V^*)\,.
\end{eqnarray*}
This gives a $GL(V)$ equivariant
decomposition of
$$\mathfrak{R}(V)=\mathfrak{U}(V)\oplus\Lambda^2(V^*)\oplus S^2(V^*)$$
as the direct sum of irreducible
$GL(V)$ modules. We have
\begin{eqnarray*}
&&\dim\{\mathfrak{U}(V)\}=\textstyle\frac13m^2(m^2-4),\quad
  \dim\{\Lambda^2(V^*)\}=\frac12m(m-1),\\
&&\dim\{S^2(V^*)\}=\textstyle\frac12m(m+1),\quad\dim\{\mathfrak{R}(V)\}=\textstyle\frac13m^2(m^2-1)\,.
\end{eqnarray*}
\end{theorem}

\begin{proof} We check the splitting as follows. If $\omega\in \Lambda^2(V^*)$,
let $\mathcal{R}_\omega:=\sigma_{\pi_a\circ\rho_{14}}(\omega)$. Then $\mathcal{R}_\omega(x,y)=-\mathcal{R}_\omega(y,x)$. 
We check the
Bianchi identity by computing:
\medbreak\qquad
$\mathcal{R}_\omega(x,y)z+\mathcal{R}_\omega(y,z)x+\mathcal{R}_\omega(z,x)y
=\textstyle\frac{-1}{1+m}\{2\omega(x,y)z+\omega(x,z)y-\omega(y,z)x$
\smallbreak\qquad\quad
$+2\omega(y,z)x+\omega(y,x)z-\omega(z,x)y$
$+2\omega(z,x)y+\omega(z,y)x-\omega(x,y)z\}$
\smallbreak\qquad\quad$=0$.
\medbreak\noindent Thus $\mathcal{R}_\omega\in\mathfrak{R}(V)$. One also has that:
\medbreak\qquad
$\rho_{14}(\mathcal{R}_\omega)(y,z)
=\textstyle\frac{-1}{1+m}\textstyle\sum_ie^i\{2\omega(e_i,y)z+\omega(e_i,z)y-\omega(y,z)e_i\}$
\smallbreak\qquad\qquad\qquad\phantom{a}
$=\textstyle\frac{-1}{1+m}\{2\omega(z,y)+\omega(y,z)-m\omega(y,z)\}=\omega(y,z)$.
\medbreak Let $\psi\in S^2(V^*)$ and let $\mathcal{R}_\psi=\sigma_{\pi_s\circ\rho_{14}}(\psi)$. Again,
$\mathcal{R}_\psi(x,y)=-\mathcal{R}_\psi(y,x)$. We verify the Bianchi identity by computing:
\medbreak\qquad
$\mathcal{R}_\psi(x,y)z+\mathcal{R}_\psi(y,z)x+\mathcal{R}_\psi(z,x)y$
\smallbreak\qquad\quad$=\frac1{1-m}\{\psi(x,z)y-\psi(y,z)x+\psi(y,x)z-\psi(z,x)y+\psi(z,y)x-\psi(x,y)z\}$
\smallbreak\qquad\quad$=0$.
\medbreak\noindent This shows that $\mathcal{R}_\psi\in\mathfrak{R}(V)$. Furthermore:
\medbreak\qquad
$\rho_{14}(\mathcal{R}_\psi)(y,z)=\textstyle\frac1{1-m}\sum_ie^i\{\psi(e_i,z)y-\psi(y,z)e_i\}$
\smallbreak\qquad\qquad\qquad\phantom{a.}$=\textstyle\frac1{1-m}\{\psi(y,z)-m\psi(y,z)\}=\psi(y,z)$.
\medbreak\noindent Consequently one has an equivariant decomposition of $\mathfrak{R}(V)$ into $GL(V)$ modules:
$$
\mathfrak{R}(V)=\mathfrak{U}(V)\oplus\Lambda^2(V^*)\oplus S^2(V^*)\,.
$$ 
We refer to \cite{S88} for the proof of the remaining
assertions of the Theorem.
\end{proof}

We say that two torsion free connections $\nabla$ and $\bar\nabla$ on a differentiable manifold $M$ are
{\it projectively equivalent} if and only if every every geodesic for $\nabla$ can be reparametrized
to be a geodesic for $\bar\nabla$, or equivalently if there exists a smooth
$1$-form $\omega$ so
$$\nabla_xy-\bar\nabla_xy=\omega(x)y+\omega(y)x\,.$$
The summand $\mathfrak{U}(V)$ plays the role of the Weyl projective tensor; it also plays a role in the
affine setting as we shall see presently in Theorem \ref{thm-4.1}. Let $\pi_{\mathfrak{U}}$ be the
associated projection on this summand in the decomposition of Theorem \ref{thm-2.3}. One has
\cite{SSV91,S88,W21}:

\begin{lemma}\label{lem-2.4}
Let $\nabla$ and $\bar\nabla$ be torsion free connections on $M$. 
\begin{enumerate}
\item If $\nabla$ and $\bar\nabla$ are projectively equivalent, then
$\pi_{\mathfrak{U}}\mathcal{R}=\pi_{\mathfrak{U}}\bar{\mathcal{R}}$.
\item The connection $\nabla$ is projectively flat if and only if $\pi_{\mathfrak{U}}\mathcal{R}=0$.
\end{enumerate}\end{lemma}

\subsection{The action of the orthogonal group on $\mathfrak{R}(V)$}\label{sect-2.4}

The associated orthogonal group $O(V,\langle\cdot,\cdot\rangle)$ acts on
$\mathfrak{R}(V)$ and on $\mathfrak{r}(V)$; the natural map from
$\mathfrak{R}(V)$ to $\mathfrak{r}(V)$ given by Equation (\ref{eqn-1.m})
is an equivariant isomorphism. Let $\Xi$ be as in Equation (\ref{eqn-1.r}). We define:
\begin{eqnarray*}
&&\Lambda^2(V^*):=\{\omega\in\otimes^2V^*:\omega_{ij}=-\omega_{ji}\},\\
&&S_0^2(V^*,\langle\cdot,\cdot\rangle):=\{\psi\in\otimes^2V^*:\psi_{ij}=\psi_{ji},\quad\textstyle\sum_{ij}\Xi^{ij}\psi_{ij}=0\},\\
&&\mathfrak{w}(V,\langle\cdot,\cdot\rangle):=\{\Theta\in\otimes^4V^*:\Theta_{ijkl}+\Theta_{jkil}+\Theta_{kijl}=0,\\
&&\qquad\qquad\qquad\qquad\textstyle\Theta_{ijkl}=-\Theta_{jikl}=\Theta_{klij},\sum_{il}\Xi^{il}\Theta_{ijkl}=0\},\\
&&\Lambda^2(\Lambda^2(V^*)):=\{\Theta\in\otimes^4V^*:\Theta_{ijkl}=-\Theta_{jikl}=-\Theta_{ijlk}=-\Theta_{klij}\},\\
&&\Lambda_0^2(\Lambda^2(V^*)):=\{\Theta\in\Lambda^2(\Lambda^2(V^*)):
\textstyle\sum_{il}\Xi^{il}\Theta_{ijkl}=0\},\\
&&\mathfrak{S}(V,\langle\cdot,\cdot\rangle):=\{\Theta\in\otimes^4V^*:\Theta_{ijkl}=-\Theta_{jikl}=\Theta_{ijlk},
  \textstyle\sum_{il}\Xi^{il}\Theta_{ijkl}=0,\\
&&\qquad\qquad\qquad\qquad\Theta_{kjil}+\Theta_{ikjl}-\Theta_{ljik}-\Theta_{iljk}=0\}\,.
\end{eqnarray*}
Note that $\Lambda^2(\Lambda^2(V^*))$, $\Lambda_0^2(\Lambda^2(V^*))$, and $\mathfrak{S}(V,\langle\cdot,\cdot\rangle)$ are not subsets
of
$\mathfrak{a}(V)$.

\begin{theorem}\label{thm-2.5} 
\ \begin{enumerate}
\item There is an
$O(V,\langle\cdot,\cdot\rangle)$ equivariant orthogonal decomposition of
$$\mathcal{R}(V)\approx\mathfrak{r}(V)=W_1\oplus...\oplus W_8$$
as the direct sum of irreducible $O(V,\langle\cdot,\cdot\rangle)$ modules where:
$$\begin{array}{ll}
\dim\{W_1\}=1,&
\dim\{W_2\}=\dim\{W_5\}=\textstyle\frac{(m-1)(m+2)}2,\\
\dim\{W_3\}=\dim\{W_4\}=\textstyle\frac{m(m-1)}2,&
\dim\{W_6\}=\textstyle\frac{m(m+1)(m-3)(m+2)}{12},\\
\dim\{W_7\}=\textstyle\frac{(m-1)(m-2)(m+1)(m+4)}8,&
\dim\{W_8\}=\textstyle\frac{m(m-1)(m-3)(m+2)}8\,.
\end{array}$$
\item There are the following isomorphisms as $O(\langle\cdot,\cdot\rangle)$
modules:
\begin{enumerate}
\item $W_1\approx\mathbb{R}$,
  $W_2\approx W_5\approx S_0^2(V^*,\langle\cdot,\cdot\rangle)$, and
  $W_3\approx W_4\approx\Lambda^2(V^*)$.
\item $W_6\approx\mathfrak{w}(V,\langle\cdot,\cdot\rangle)$ is the space of Weyl conformal curvature tensors.
\item  $W_7\approx\mathfrak{S}(V,\langle\cdot,\cdot\rangle)$ and $W_8\approx\Lambda_0^2(\Lambda^2(V^*))$.
\end{enumerate}
\end{enumerate}
\end{theorem}

We refer to Bokan \cite{B90} for the proof of Assertion (1) in
the context of a positive definite inner product; it extends immediately to the
indefinite inner products. We will prove Assertion (2a) later in this section. We will prove Assertion (2b) in Section
\ref{sect-3}. We will prove Assertion (2c) in Section \ref{sect-4}.

\begin{remark}\label{thm-2.6}
\rm  Since $W_2$ and $W_5$ are isomorphic as $O(V,\langle\cdot,\cdot\rangle)$
modules and since $W_3$ and $W_4$ are isomorphic as $O(V,\langle\cdot,\cdot\rangle)$ modules, the decomposition of
$\mathfrak{R}(V)$ into irreducible module summands is not unique; this fact plays an important role in the analysis of Bokan
\cite{B90}.\end{remark}

We shall need a technical result before proving Theorem \ref{thm-2.5} (2). We use Equation (\ref{eqn-1.m}) to lower indices
and to define a curvature tensor $R$ associated to a curvature operator $\mathcal{R}$. Let $\Xi$ be as in Equation
(\ref{eqn-1.r}). Then:
\begin{eqnarray*}
&&\rho_{14}(R)(x,y):=\textstyle\sum_{ij}\Xi^{ij}R(e_i,x,y,e_j),\quad
  \rho_{23}(R)(x,y):=\textstyle\sum_{ij}\Xi^{ij}R(x,e_i,e_j,y),\\
&&\rho_{24}(R)(x,y):=\textstyle\sum_{ij}\Xi^{ij}R(x,e_i,y,e_j),\quad
  \rho_{13}(R)(x,y):=\textstyle\sum_{ij}\Xi^{ij}R(e_i,x,e_j,y),\\
&&\rho_{34}(R)(x,y):=\textstyle\sum_{ij}\Xi^{ij}R(x,y,e_i,e_j)=-\rho_{14}(R)(x,y)+\rho_{14}(R)(y,x)\,.
\end{eqnarray*}

There is an $O(V,\langle\cdot,\cdot\rangle)$ equivariant decomposition:
$$V^*\otimes V^*=\Lambda^2(V^*)\oplus S_0^2(V^*,\langle\cdot,\cdot\rangle)
\oplus\mathbb{R}$$
where $S_0^2(V^*,\langle\cdot,\cdot\rangle)$ is the space of trace free symmetric bilinear forms, and where
$\mathbb{R}$ is the trivial $O(V,\langle\cdot,\cdot\rangle)$ module. If
$\pi_a$, $\pi_0$, and $\tau$ are the associated orthogonal projections,
then
\begin{equation}\label{eqn-2.c}
\begin{array}{l}
\pi_a(\psi)(x,y):=\textstyle\frac12\{\psi(x,y)-\psi(y,x)\},\\
\pi_s(\psi)(x,y):=\textstyle\frac12\{\psi(x,y)+\psi(y,x)\},\vphantom{\vrule height 11pt}\\
\tau(\psi):=\textstyle\sum_{ij}\Xi^{ij}\psi(e_i,e_j),\vphantom{\vrule height 11pt}\\
\pi_0(\psi)(x,y):=\pi_s(\psi)(x,y)
   -\frac1m\tau(\psi)\langle\cdot,\cdot\rangle\,.\vphantom{\vrule height 11pt}
\end{array}\end{equation}
There is only one non-trivial scalar curvature arising from $R\in\mathfrak{r}(V)$ since
$$
\begin{array}{l}
\tau(\rho_{14}(R))=\textstyle\sum_{ijkl}\Xi^{il}\Xi^{jk}R(e_i,e_j,e_k,e_l)=
  \tau(\rho_{23}(R))=-\tau(\rho_{24}(R)),\\
\tau(\rho_{34}(R))=\textstyle\sum_{ijkl}\Xi^{ij}\Xi^{kl}R(e_i,e_j,e_k,e_l)=0\,.
\end{array}
$$

If $\psi\in S_0^2(V^*,\langle\cdot,\cdot\rangle)$ and if $\omega\in\Lambda^2(V^*)$, let:
\begin{eqnarray*}
&&\sigma_1(\psi)(x,y,z,w):=\psi(x,w)\langle y,z\rangle-\psi(y,w)\langle x,z\rangle,\\
&&\sigma_2(\psi)(x,y,z,w):=\langle x,w\rangle\psi(y,z)-\langle y,w\rangle\psi(x,z),\\
&&\sigma_3(\omega)(x,y,z,w):=2\omega(x,y)\langle z,w\rangle+\omega(x,z)\langle y,w\rangle-\omega(y,z)\langle x,w\rangle,\\
&&\sigma_4(\omega)(x,y,z,w):=\omega(x,w)\langle y,z\rangle-\omega(y,w)\langle x,z\rangle\,.
\end{eqnarray*}

\begin{lemma}\label{lem-2.7}
\ \begin{enumerate}
\item $\sigma_1$ and $\sigma_2$ are $O(V,\langle\cdot,\cdot\rangle)$ equivariant maps from $S_0^2(V^*,\langle\cdot,\cdot\rangle)$
to
$\mathfrak{r}(V)$, $\sigma_3$ and $\sigma_4$ are $O(V,\langle\cdot,\cdot\rangle)$ equivariant
maps from $\Lambda^2(V^*)$ to $\mathfrak{r}(V)$, and\begin{eqnarray*}
&&\left(\begin{array}{ll}
\rho_{14}\circ\sigma_1&\rho_{23}\circ\sigma_1\\
\rho_{14}\circ\sigma_2&\rho_{23}\circ\sigma_2\end{array}\right)
=\left(\begin{array}{rr}
-\id&(m-1)\id\\(m-1)\id&-\id\end{array}\right),\\
&&\left(\begin{array}{rr}
\rho_{13}\circ\sigma_3&\rho_{34}\circ\sigma_3\\
\rho_{13}\circ\sigma_4&\rho_{34}\circ\sigma_4
\end{array}\right)
=\left(\begin{array}{rr}-3\id&2(m+1)\id\\(1-m)\id&2\id\end{array}\right)\,.
\end{eqnarray*}
\item We have $O(V,\langle\cdot,\cdot\rangle)$ equivariant sequences which are equivariantly split:
\begin{eqnarray*}
&&\tau\circ\rho_{14}:\mathfrak{r}(V)\rightarrow\mathbb{R}\rightarrow0,\\
&&\pi_0\circ\rho_{14}\oplus\pi_0\circ\rho_{13}:\mathfrak{r}(V)\rightarrow
S_0^2(V^*,\langle\cdot,\cdot\rangle)\oplus S_0^2(V^*,\langle\cdot,\cdot\rangle)\rightarrow0,\\
&&\pi_a\circ\rho_{13}\oplus\pi_a\circ\rho_{34}:\mathfrak{r}(V)\rightarrow
\Lambda^2(V^*)\oplus\Lambda^2(V^*)\rightarrow0\,.
\end{eqnarray*}
\end{enumerate}
\end{lemma}

\begin{proof}
Let $\psi\in S_0^2(V^*,\langle\cdot,\cdot\rangle)$ and let $\omega\in\Lambda^2(V^*)$. Set $R_1:=\sigma_1(\psi)$, $R_2:=\sigma_2(\psi)$,
$R_3:=\sigma_3(\omega)$, and
$R_4:=\sigma_4(\omega)$. It is immediate $R_i(x,y,z,w)=-R_i(y,x,z,w)$. To show that $R_i\in\mathfrak{r}(V)$, we must verify the
first Bianchi identity is satisfied:
\medbreak\qquad $R_1(x,y,z,w)+R_1(y,z,x,w)+R_1(z,x,y,w)$
\par\quad\qquad $=\psi(x,w)\langle y,z\rangle-\psi(y,w)\langle x,z\rangle$\par\qquad\quad
\phantom{.}$+\psi(y,w)\langle z,x\rangle-\psi(z,w)\langle y,x\rangle$\par\qquad\quad
\phantom{.}$+\psi(z,w)\langle x,y\rangle-\psi(x,w)\langle z,y\rangle=0$,
\medbreak\qquad
$R_2(x,y,z,w)+R_2(y,z,x,w)+R_2(z,x,y,w)$
\par\qquad\quad$=\langle x,w\rangle\psi(y,z)-\langle y,w\rangle\psi(x,z)$
\par\qquad\quad\phantom{.}$ +\langle y,w\rangle\psi(z,x)-\langle z,w\rangle\psi(y,x)$
\par\qquad\quad\phantom{.}$+\langle z,w\rangle\psi(x,y)-\langle x,w\rangle\psi(z,y)=0$,
\medbreak\qquad
$R_3(x,y,z,w)+R_3(y,z,x,w)+R_3(z,x,y,w)$
\par\qquad\quad$=2\omega(x,y)\langle z,w\rangle
+\omega(x,z)\langle y,w\rangle-\omega(y,z)\langle x,w\rangle$
\par\qquad\quad\phantom{.}$+2\omega(y,z)\langle x,w\rangle
+\omega(y,x)\langle z,w\rangle-\omega(z,x)\langle y,w\rangle$
\par\qquad\quad\phantom{.}$+2\omega(z,x)\langle y,w\rangle
+\omega(z,y)\langle x,w\rangle-\omega(x,y)\langle z,w\rangle=0$,
\medbreak\qquad$R_4(x,y,z,w)+R_4(y,z,x,w)+R_4(z,x,y,w)$
\par\qquad\quad$=\omega(x,w)\langle y,z\rangle-\omega(y,w)\langle x,z\rangle$
\par\qquad\quad\phantom{.}$+\omega(y,w)\langle z,x\rangle-\omega(z,w)\langle y,x\rangle$
\par\qquad\quad\phantom{.}$+\omega(z,w)\langle x,y\rangle-\omega(x,w)\langle z,y\rangle=0$.

\medbreak We complete the proof of Assertion (1) by computing:
\medbreak\quad
$\rho_{14}(R_1)(y,z)=\textstyle\sum_{ij}\Xi^{ij}\{\psi(e_i,e_j)\langle y,z\rangle-\psi(y,e_j)\langle e_i,z\rangle\}$
\par\qquad\quad
$=\tau(\psi)\langle y,z\rangle-\psi(y,z)=-\psi(y,z)$,
\medbreak\quad
$\rho_{23}(R_1)(x,w)=\textstyle\sum_{ij}\Xi^{ij}\{\psi(x,w)\langle e_i,e_j\rangle-\psi(e_i,w)\langle x,e_j\rangle\}$
\par\qquad\quad
$=(m-1)\psi(x,w)$,
\medbreak\quad
$\rho_{14}(R_2)\psi(y,z)=\textstyle\sum_{ij}\Xi^{ij}\{\langle e_i,e_j\rangle\psi(y,z)-\langle y,e_j\rangle\psi(e_i,z)\}$
\par\qquad\quad$=(m-1)\psi(y,z)$,
\medbreak\quad
$\rho_{23}(R_2)(x,w)=\textstyle\sum_{ij}\Xi^{ij}\{\langle x,w\rangle\psi(e_i,e_j)-\langle e_i,w\rangle\psi(x,e_j)\}$
\par\qquad\quad$=\tau(\psi)\langle x,w\rangle-\psi(x,w)=-\psi(x,w)$,
\medbreak\quad
$\rho_{13}(R_3)(y,w)=\textstyle\sum_{ij}\Xi^{ij}\{2\omega(e_i,y)\langle e_j,w\rangle+\omega(e_i,e_j)\langle
y,w\rangle-\omega(y,e_j)\langle e_i,w\rangle$
\par\qquad\quad$=-3\omega(y,w)$,
\medbreak\quad
$\rho_{34}(R_3)(x,y)=\textstyle\sum_{ij}\Xi^{ij}\{2\omega(x,y)\langle e_i,e_j\rangle+\omega(x,e_i)\langle
y,e_j\rangle-\omega(y,e_i)\langle x,e_i\rangle$
\par\qquad\quad$=2(m+1)\omega(x,y)$,
\medbreak\quad
$\rho_{13}(R_4)(y,w)=\textstyle\sum_{ij}\Xi^{ij}\{\omega(e_i,w)\langle y,e_j\rangle-\omega(y,w)\langle e_i,e_j\rangle\}$
$=(1-m)\omega(y,w)$,
\medbreak\quad
$\rho_{34}(R_4)(x,y)=\textstyle\sum_{ij}\Xi^{ij}\{\omega(x,e_j)\langle y,e_i\rangle-\omega(y,e_j)\langle x,e_i\rangle\}$
$=2\omega(x,y)$.

\medbreak We now prove Assertion (2). We show the first sequence splits by computing:
\begin{eqnarray*}
\textstyle\frac1{m(m-1)}\tau(\rho_{14}(\sigma_1\langle\cdot,\cdot\rangle))&=&
\textstyle\frac1{m(m-1)}\textstyle\sum_{ijkl}\Xi^{il}\Xi^{jk}\{\Xi_{il}\Xi_{jk}-\Xi_{ik}\Xi_{jl}\}\\
&=&\textstyle\frac1{m(m-1)}\sum_{ij}\{\delta_i^i\delta_j^j-\delta_i^j\delta_j^i\}=1\,.
\end{eqnarray*}
As the determinants of the two coefficient matrices in Assertion (1) are non-zero, the desired
splitting of the second and of the third sequences follows.
\end{proof}

\medbreak\noindent{\it Proof of Theorem \ref{thm-2.5} (2a)}. By Lemma \ref{lem-2.7}, $\mathbb{R}$ has multiplicity 1,
$S_0^2(V^*,\langle\cdot,\cdot\rangle)$ has multiplicity $2$, and
$\Lambda^2(V^*)$ has multiplicity $2$ in the decomposition of $\mathfrak{r}(V)$ as an $O(\langle\cdot,\cdot\rangle)$
module. These modules are irreducible and
$$\dim\{\mathbb{R}\}=1,\quad
  \dim\{S_0^2(V^*,\langle\cdot,\cdot\rangle)\}=\textstyle\frac{(m-1)(m+2)}2,\quad
  \dim\{\Lambda^2(V^*)\}=\textstyle\frac{m(m-1)}2\,.$$
Theorem \ref{thm-2.5} (2a) now follows from Theorem \ref{thm-2.5} (1).
\hfill\qedbox

\section{Algebraic Curvature Tensors}\label{sect-3}
In this section, we study the quadralinear forms with the
same symmetries as those of the Levi-Civita connection of a pseudo-Riemannian manifold.

\subsection{The action of the general linear group on $\mathfrak{a}(V)$}\label{sect-3.1}

\begin{theorem}\label{thm-3.1} 
$\mathfrak{a}(V)$ is an irreducible $GL(V)$ module.
\end{theorem}

We postpone the proof of this result until Section \ref{sect-5} as we must first establish some additional notation. 
\subsection{The action of $O(V,\langle\cdot,\cdot\rangle)$ on $\mathfrak{a}(V)$}\label{sect-3.2}
Let
\begin{equation}\label{eqn-3.a}
(\id\otimes\pi_s)(R)(x,y,z,w):=\textstyle\frac12\{R(x,y,z,w)+R(x,y,w,z)\}\text{ for }R\in\mathfrak{r}(V)\,.
\end{equation}
If $\phi,\psi\in S^2(V^*)$, one can define an algebraic
curvature tensor
$\phi\wedge\psi\in\mathfrak{a}(V)$ by:
\begin{equation}\label{eqn-3.b}
\begin{array}{ll}
\{\phi\wedge\psi\}(x,y,z,w):=&\textstyle\frac12\{
  \phi(x,w)\psi(y,z)-\phi(x,z)\psi(y,w)\\
&\phantom{}+\phi(y,z)\psi(x,w)-\phi(y,w)\psi(x,z)\}\,.
\end{array}\end{equation}
(This has a different normalizing constant than the usual Kulkarni-Nomizu product). These tensors arise naturally. If $L$ is the
second fundamental form of a hypersurface
$M$ in
$\mathbb{R}^{m+1}$, then
$$R_M=L\wedge L\,.$$
Define:
\begin{equation}\label{eqn-3.c}
\begin{array}{l}
\mathfrak{w}(V,\langle\cdot,\cdot\rangle):=\ker\{\rho_{14}\}\cap\mathfrak{a}(V),\\
\sigma_{\id\otimes\pi_s}(S)_{ijkl}:=S_{ijkl}+\textstyle\frac12\{S_{kjil}+S_{ikjl}-S_{ljik}-S_{iljk}\},\\
\vphantom{\vrule height 11pt}
\sigma_{\mathfrak{a},\rho_{14}}(\psi):=\textstyle\frac2{m-2}\psi\wedge\langle\cdot,\cdot\rangle-\frac{\tau(\psi)}{(m-1)(m-2)}
 \langle\cdot,\cdot\rangle\wedge\langle\cdot,\cdot\rangle\,.\vphantom{\vrule height 11pt}
\end{array}\end{equation}

\begin{theorem}\label{thm-3.2}
\ \begin{enumerate}
\item Let $R\in\otimes^4V^*$ satisfy Equations
(\ref{eqn-1.c}) and (\ref{eqn-1.d}). Then Equations (\ref{eqn-1.e}) and (\ref{eqn-1.f}) are
equivalent.
\item The maps $\id\otimes\pi_s$ and $\rho_{14}$ define $GL(V)$ and $O(V,\langle\cdot,\cdot\rangle)$ equivariant short exact
sequences, respectively,
\begin{eqnarray*}
&&0\rightarrow\mathfrak{a}(V)\rightarrow\mathfrak{r}(V)\mapright{\id\otimes\pi_s}\Lambda^2(V^*)\otimes S^2(V^*)\rightarrow0,\\
&&0\rightarrow\mathfrak{w}(V,\langle\cdot,\cdot\rangle)\rightarrow\mathfrak{a}(V)\mapright{\rho_{14}}S^2(V^*)\rightarrow0\,.
\end{eqnarray*}
which are equivariantly split, respectively, by the maps $\sigma_{\id\otimes\pi_s}$ and $\sigma_{\mathfrak{a},\rho_{14}}$.
\item This gives an $O(V,\langle\cdot,\cdot\rangle)$ equivariant decomposition of
$$\mathfrak{a}(V)\approx\mathfrak{w}(V,\langle\cdot,\cdot\rangle)\ \oplus\ 
  S_0^2(V^*,\langle\cdot,\cdot\rangle)\ \oplus\ \{\mathbb{R}\}$$
as the direct sum of irreducible $O(V,\langle\cdot,\cdot\rangle)$ modules where
$$\begin{array}{ll}
\dim\{\mathfrak{w}(V,\langle\cdot,\cdot\rangle)\}=\textstyle\frac1{12}m(m+1)(m+2)(m-3),&\dim\{\mathbb{R}\}=1,\\
\dim\{S_0^2(V^*,\langle\cdot,\cdot\rangle)\}=\frac12(m-1)(m+2),&\dim\{\mathfrak{a}(V)\}=\textstyle\frac1{12}m^2(m^2-1)\,.
\end{array}$$
\end{enumerate}
\end{theorem}

\begin{proof} It is immediate that (\ref{eqn-1.c}) and (\ref{eqn-1.e}) imply Equation
(\ref{eqn-1.f}). Conversely, suppose that Equations (\ref{eqn-1.c}), (\ref{eqn-1.d}), and
(\ref{eqn-1.f}) hold. We use the following notation:
\begin{eqnarray*}
&&R(\xi_1,\xi_2,\xi_3,\xi_4)=a_1,\quad R(\xi_3,\xi_4,\xi_1,\xi_2)=a_1+\varepsilon_1,\\
&&R(\xi_1,\xi_3,\xi_2,\xi_4)=a_2,\quad R(\xi_2,\xi_4,\xi_1,\xi_3)=a_2+\varepsilon_2,\\
&&R(\xi_2,\xi_3,\xi_1,\xi_4)=a_3,\quad R(\xi_1,\xi_4,\xi_2,\xi_3)=a_3+\varepsilon_3\,.
\end{eqnarray*}
We establish Assertion (1) by showing $\varepsilon_1=\varepsilon_2=\varepsilon_3=0$. We compute:
\medbreak\qquad 
$0=R(\xi_1,\xi_2,\xi_3,\xi_4)+R(\xi_2,\xi_3,\xi_1,\xi_4)+R(\xi_3,\xi_1,\xi_2,\xi_4)$
\par\qquad\quad$=a_1+a_3-a_2$,
\medbreak\qquad $0=R(\xi_1,\xi_2,\xi_4,\xi_3)+R(\xi_2,\xi_4,\xi_1,\xi_3)+R(\xi_4,\xi_1,\xi_2,\xi_3)$
\par\qquad\quad$=-a_1+a_2-a_3+\varepsilon_2-\varepsilon_3=\varepsilon_2-\varepsilon_3$,
\medbreak\qquad$0=R(\xi_1,\xi_3,\xi_4,\xi_2)+R(\xi_3,\xi_4,\xi_1,\xi_2)+R(\xi_4,\xi_1,\xi_3,\xi_2)$
\par\qquad\quad$=-a_2+a_1+a_3+\varepsilon_1+\varepsilon_3=\varepsilon_1+\varepsilon_3$,
\medbreak\qquad$0=R(\xi_2,\xi_3,\xi_4,\xi_1)+R(\xi_3,\xi_4,\xi_2,\xi_1)+R(\xi_4,\xi_2,\xi_3,\xi_1)$
\par\qquad\quad$=-a_3-a_1+a_2-\varepsilon_1+\varepsilon_2
  =-\varepsilon_1+\varepsilon_2$.
\medbreak This yields the equations
$0=\varepsilon_2-\varepsilon_3=\varepsilon_1+\varepsilon_3=
  -\varepsilon_1+\varepsilon_2
$
from which it follows that
$\varepsilon_1=\varepsilon_2=\varepsilon_3=0$; this proves Assertion (1).

 Let $S\in\Lambda^2(V^*)\otimes
S^2(V^*)$. We compute:
\medbreak
\qquad$\sigma_{\id\otimes\pi_s}(S)_{ijkl}+\sigma_{\id\otimes\pi_s}(S)_{jikl}$
\smallbreak\qquad\qquad$=S_{ijkl}+\textstyle\frac12(S_{kjil}+S_{ikjl}-S_{ljik}-S_{iljk})$
\smallbreak\qquad\qquad\phantom{.}$+S_{jikl}+\frac12(S_{kijl}+S_{jkil}-S_{lijk}-S_{jlik})=0$,
\medbreak\qquad
$\sigma_{\id\otimes\pi_s}(S)_{ijkl}+\sigma_{\id\otimes\pi_s}(S)_{jkil}+\sigma_{\id\otimes\pi_s}(S)_{kijl}$
\smallbreak\qquad\qquad$=S_{ijkl}+\textstyle\frac12(S_{kjil}+S_{ikjl}-S_{ljik}-S_{iljk})$
\smallbreak\qquad\qquad\phantom{.}$+S_{jkil}+\textstyle\frac12(S_{ikjl}+S_{jikl}-S_{lkji}-S_{jlki})$
\smallbreak\qquad\qquad\phantom{.}$+S_{kijl}+\textstyle\frac12(S_{jikl}+S_{kjil}-S_{likj}-S_{klij})$
\smallbreak\qquad\qquad$=0$.
\medbreak\noindent This shows that $\sigma_{\id\otimes\pi_s}$ takes values in $\mathfrak{r}(V)$. 
Let $\alpha(S):=\sigma_{\id\otimes\pi_s}S-S$. Then
\begin{equation}\label{eqn-3.d}
\alpha(S)_{ijkl}:=\textstyle\frac12(S_{kjil}+S_{ikjl}-S_{ljik}-S_{iljk})\in\Lambda^2(V^*)\otimes\Lambda^2(V^*)\,.
\end{equation}
The map $\alpha$ will also play a role in Section \ref{sect-4.2}. Since $\id\otimes\pi_s$ vanishes
on the space $\Lambda^2(V^*)\otimes\Lambda^2(V^*)$, one has that
$$(\id\otimes\pi_s)(\sigma_{\id\otimes\pi_s}(S))=(\id\otimes\pi_s)(S)+(\id\otimes\pi_s)\alpha(S)=S\,.$$
This shows that $\id\otimes\pi_s$  is an equivariant splitting.
We refer to Singer and Thorpe \cite{ST} or to Strichartz \cite{S88} for the proof of
the remaining assertions.
\end{proof}

\medbreak\noindent{\it Proof of Theorem \ref{thm-2.5} (2b)}. Because $\mathfrak{w}$ is the space of {\it Weyl conformal
tensors},
$$\dim\{\mathfrak{w}(V,\langle\cdot,\cdot\rangle)\}=\textstyle\frac1{12}m(m+1)(m-3)(m+2)=\dim\{W_6\}\,.$$
Since $\mathfrak{w}(V,\langle\cdot,\cdot\rangle$ is an irreducible $O(V,\langle\cdot,\cdot\rangle)$ module,
we may use Theorem \ref{thm-2.5} (1) to identify $W_6=\mathfrak{w}(V,\langle\cdot,\cdot\rangle)$.
\hfill\qedbox

\medbreak Theorem \ref{thm-2.1} generalizes immediately to this setting:

\begin{theorem}\label{thm-3.3}
Let $A\in\mathfrak{a}(V)$ be given. Regard $V$ as a smooth manifold in its own right. Let $0$ be the origin
of $V$ and identify $T_0V=V$. There exists a pseudo-Riemannian metric $g$ defined on $V$ so that
$R_0^g=A$ where $R_0^g$ is the curvature tensor of the associated Levi-Civita connection.
\end{theorem}

\begin{proof} Let
$\{e_i\}$ be an orthonormal basis for $V$. Let $x_i$ be the associated coordinate system. We define the germ
of a pseudo-Riemannian metric on $V$ by setting
$$
g_{ab}=g(\partial_{x_a},\partial_{x_b}):=\langle e_a,e_b\rangle-\textstyle\frac13\textstyle\sum_{cd}A_{acdb}x_cx_d\,.
$$
Clearly $g_{ab}=g_{ba}$. As $g|_{T_0V}=\langle\cdot,\cdot\rangle$, $g$ is non-degenerate near $0$. One may then use a partition
of unity to extend $g$ to be non-degenerate on all of $V$ without changing it near $0$. One has
$$\Gamma_{ijk}:=g(\nabla_{\partial_{x_i}}\partial_{x_j},\partial_{x_k})=
   \ffrac12(\partial_{x_i}g_{jk}+\partial_{x_j}g_{ik}-\partial_{x_k}g_{ij})\,.
$$
Since $\Gamma_{ijk}(0)=0$, one has
\begin{eqnarray*}
&&R_{ijkl}(0):=R(\partial_{x_i},\partial_{x_j},\partial_{x_k},\partial_{x_l})(0)
 =\{\partial_{x_i}\Gamma_{jkl}-\partial_{x_j}\Gamma_{ikl}\}(0)\\
&=&\textstyle\frac12\{\partial_{x_i}(\partial_{x_j}g_{kl}+\partial_{x_k}g_{jl}-\partial_{x_l}g_{jk})
   -\partial_{x_j}(\partial_{x_i}g_{kl}+\partial_{x_k}g_{il}-\partial_{x_l}g_{ik})\}(0)
\\&=&
 \textstyle\frac16\{-A_{jikl}-A_{jkil}+A_{jilk}+A_{jlik}+A_{ijkl}+A_{ikjl}-A_{ijlk}-A_{iljk}
   \}\\
  &=&\textstyle\frac16\{4A_{ijkl}-2A_{iljk}-2A_{iklj}\}
 =A_{ijkl}\,.
\end{eqnarray*}
The desired result now follows.\end{proof}

The following result was first proved by Fiedler \cite{F01} using Young symmetrizers;
subsequently Gilkey \cite{G02} established it using a direct construction and D\'{\i}az-Ramos and
Garc\'{\i}a-R\'{\i}o \cite{DRGR-04} derived it from the Nash embedding theorem. We adopt the notation of Equation (\ref{eqn-3.b})
to define $\phi\wedge\psi\in\mathfrak{a}(V)$ for $\phi,\psi\in S^2(V^*)$.

\begin{theorem}\label{thm-3.4}
$\mathfrak{a}(V)=\operatorname{Span}_{\mathbb{R}}\{\phi\wedge\phi:\phi\in
S^2(V^*)\}$.
\end{theorem}

We use Theorem \ref{thm-3.2} to establish a slightly stronger version of Theorem \ref{thm-3.4}:

\begin{theorem}\label{thm-3.5}
\ \begin{enumerate}
\item If $A\in\mathfrak{a}(V)$, there is a finite collection of elements $\phi_\nu\in S^2(V^*)$ such
that
$\Rank\{\phi_\nu\}=2$ and such that
$A=\sum_\nu\phi_\nu\wedge\phi_\nu$.
\item Suppose given $(p,q)$ with $2\le p+q\le m$. Let $S^2_{(p,q)}(V^*)$ be the set of all symmetric bilinear forms on $V$
of signature $(p,q)$. Then 
$$
\mathfrak{a}(V)=\operatorname{Span}_{\phi\in S^2_{(p,q)}(V^*)}\{\phi\wedge\phi\}\,.
$$
\end{enumerate}
\end{theorem}

\begin{proof} Consider the following
$GL(V)$ invariant subspace of
$\mathfrak{a}(V)$:
$$\mathfrak{b}(V):=\operatorname{Span}_{\mathbb{R}}\{\phi\wedge\phi:\phi\in
S^2(V^*), \Rank\{\phi\}=2\}\,.
$$ 
We apply Theorem \ref{thm-3.1} to show $\mathfrak{b}(V)=\mathfrak{a}(V)$. This shows that we may express any
$A\in\mathfrak{a}(V)$ in the form
$c_1\phi_1\wedge\phi_1+...+c_k\phi_k\wedge\phi_k$ where the $\phi_\nu$ are symmetric bilinear forms of
rank $2$ and where the $c_\nu\in\mathbb{R}$. By rescaling the $\phi_\nu$, we may assume that the
$c_\nu=\pm1$. Set
$\alpha_1:=e^1\otimes e^1+e^2\otimes e^2$ and $\alpha_2:=e^1\otimes e^2+e^2\otimes e^1$. We have
$$(\alpha_1\wedge\alpha_1)(e_1,e_2,e_2,e_1)=+1\quad\text{and}\quad
  (\alpha_2\wedge\alpha_2)(e_1,e_2,e_2,e_1)=-1\,.
$$
Thus $\alpha_1\wedge\alpha_1=-\alpha_2\wedge\alpha_2$. Consequently, by replacing a definite form by an indefinite form or an
indefinite form by a definite form if necessary, we can change the sign and assume that all the constants $c_\nu$ are equal to
$1$. Assertion (1) now follows.

To prove Assertion (2), we set
$$
\mathfrak{b}(V):=\operatorname{Span}_{\phi\in S^2_{(p,q)}(V^*)}\{\phi\wedge\phi\}\,.
$$
As this is a non-empty $GL(V)$ invariant subspace of $\mathfrak{a}(V)$, Theorem \ref{thm-3.1} shows
$\mathfrak{a}(V)=\mathfrak{b}(V)$ as desired.\end{proof}

\section{Affine Curvature Tensors in the Algebraic Setting}\label{sect-4}\ 

\subsection{The action of the general linear group on $\mathfrak{F}(V)$}\label{sect-4.1}
We adopt the notion of Equation (\ref{eqn-2.a}) to define $\mathfrak{U}(V)$; the geometrical significance of this
subspace is given in Lemma \ref{lem-2.4}.

Use Equations (\ref{eqn-1.g}) and (\ref{eqn-2.b}) to define $\rho_{14}$, $\pi_a$, and $\pi_s$. Let
$\sigma_{\pi_a\circ\rho_{14}}$ and $\sigma_{\pi_s\circ\rho_{14}}$ be as in Theorem \ref{thm-2.3}. The following is an
immediate consequence of Theorem
\ref{thm-2.3}:
\begin{theorem}\label{thm-4.1}
We have the following $GL(V)$ equivariant short exact sequences
\begin{eqnarray*}
&&0\rightarrow\mathfrak{F}(V)\rightarrow\mathfrak{R}(V)\mapright{\pi_a\circ\rho_{14}}\Lambda^2(V^*)\rightarrow0,\\
&&0\rightarrow\mathfrak{U}(V)\rightarrow\mathfrak{F}(V)\mapright{\pi_s\circ\rho_{14}}S^2(V^*)\rightarrow0
\end{eqnarray*}
which are equivariantly split by the maps $\sigma_{\pi_a\circ\rho_{14}}$ and $\sigma_{\pi_s\circ\rho_{14}}$, respectively.
This gives a $GL(V)$ equivariant decomposition of 
$$\mathfrak{F}(V)=\mathfrak{U}(V)\oplus S^2(V^*)$$ as the direct sum
of irreducible $GL(V)$ modules where
\begin{eqnarray*}
&&\dim\{\mathfrak{U}(V)\}=\textstyle\frac{m^2(m^2-4)}{3},\\
&&\dim\{\dim\{S^2(V^*)\}=\textstyle\frac12m(m+1),\\
&&\dim\{\mathfrak{F}(V)\}=\textstyle\frac{m(m-1)(2m^2+2m-3)}{6}\,.
\end{eqnarray*}
\end{theorem}

We use this result to generalize Theorem \ref{thm-3.4} to the setting at hand. We exploit in an essential way that the space
$\mathfrak{A}(V,\langle\cdot,\cdot\rangle)$ depends non-trivially on the particular bilinear form which is chosen. Let
$\mathcal{G}_{(p,q)}(V)$ be the set of non-degenerate bilinear forms on $V$ of signature $(p,q)$.
Let $\mathcal{G}_{(p,q)}(M)$ be the set of all pseudo-Riemannian metrics on a smooth $m$-dimensional manifold $M$ of signature
$(p,q)$. If
$g\in\mathcal{G}_{(0,m)}(M)$ and if
$P\in M$, let
$\mathcal{R}(g,P)$ be the curvature operator of the Levi-Civita connection defined by $g$.

\begin{theorem}\label{thm-4.2}
\ \begin{enumerate}\item
If $p+q=m$, then
$\mathfrak{F}(V)=\operatorname{Span}_{\langle\cdot,\cdot\rangle\in\mathcal{G}_{(p,q)}}
   \{\mathfrak{A}(V,\langle\cdot,\cdot\rangle)\}$.
\item We have that
$\mathfrak{F}(T_PM)=\operatorname{Span}_{g\in\mathcal{G}_{(0,m)}(M)}\{\mathcal{R}(g,P)\}$.
\end{enumerate}
\end{theorem}

\begin{proof} Let
$$\mathfrak{B}(V):=\operatorname{Span}_{\langle\cdot,\cdot\rangle\in\mathcal{G}_{(p,q)}}
\{\mathfrak{A}(V,\langle\cdot,\cdot\rangle)\}\,.$$
Let $\Psi\in GL(V)$. If
$\mathcal{A}\in\mathfrak{A}(V,\langle\cdot,\cdot\rangle)$, then
$$\Psi^*\mathcal{A}\in\mathfrak{A}(V,\Psi^*\langle\cdot,\cdot\rangle)\,.$$
Thus
$\mathfrak{B}(V)$ is invariant under the action of
$GL(V)$. Since $\mathfrak{B}(V)\ne\{0\}$, Theorem \ref{thm-4.1} shows exactly one of the following alternatives holds:
\begin{enumerate}
\item $\mathfrak{B}(V)=\ker\{\pi_s\circ\rho_{14}\}$.
\item $\mathfrak{B}(V)\approx S^2(V^*)$.
\item $\mathfrak{B}(V)=\mathfrak{F}(V)$.
\end{enumerate}

 If
$\langle\cdot,\cdot\rangle\in\mathcal{G}_{(p,q)}(V)$, let
$\mathcal{A}_{\langle\cdot,\cdot\rangle}\in\mathfrak{A}(V,\langle\cdot,\cdot\rangle)$ be the associated algebraic curvature
operator of constant sectional curvature:
$$\mathcal{A}_{\langle\cdot,\cdot\rangle}(x,y)z:=\langle y,z\rangle x-\langle x,z\rangle y\,.$$
Since $\rho_{14}(\mathcal{A}_{\langle\cdot,\cdot\rangle})=(m-1)\langle\cdot,\cdot\rangle$, $\mathfrak{B}(V)\ne\ker\{\rho_{14}\}$.
This eliminates the first possibility.
Since
$m\ge4$, $m(m+1)>6$. Consequently,
$$\dim\{\mathfrak{B}(V)\}\ge\dim\{\mathfrak{A}(V,\langle\cdot,\cdot\rangle)\}=\textstyle\frac{m^2(m^2-1)}{12}
>\frac{m(m+1)}2=\dim\{S^2(V^*)\}\,.$$
This eliminates the second possibility. Thus the third possibility holds; this proves Assertion (1).

Let $V=T_PM$. Let $g_0\in\mathcal{G}_{(0,m)}(T_PM)$. By Theorem \ref{thm-3.3},
$$
\mathfrak{A}(V,g_0)=\cup_{g\in\mathcal{G}_{(0,m)},g|_{T_PM}=g_0}\{\mathcal{R}(g,P)\}\,.
$$
Assertion (2) now follows from
Assertion (1).
\end{proof}

\subsection{Centro affine geometry} Let $h\in S^2(V^*)$ and let $\mathcal{C}\in S^2(V^*)\otimes V$. Define
\begin{eqnarray*}
&&\mathcal{R}_h(x,y)z:=h(y,z)x-h(x,z)y,\\
&&\mathcal{R}_{\mathcal{C}}(w,v)u:=\mathcal{C}(v,\mathcal{C}(w,u))-\mathcal{C}(w,\mathcal{C}(v,u))\,.
\end{eqnarray*}
The decomposition of Theorem \ref{thm-4.1} has geometric
significance. Let $h$ be the centroaffine metric, let $\nabla$ be the induced connection, and
let $\nabla^*$ be the conormal connection. Then $\mathcal{R}_h$ is the curvature operator of
both $\nabla$ and of
$\nabla^*$ while the Riemannian curvature
tensor of the associated Levi-Civita connection is given
by $\mathcal{R}_C+\mathcal{R}_h$.

\begin{theorem}\label{thm-4.3}
\ 
\begin{enumerate}
\item $\mathcal{R}_h\in\sigma_{\pi_s\circ\rho_{14}}S^2(V^*)$ and $\sigma_{\pi_s\circ\rho_{14}}S^2(V^*)=\operatorname{Span}_{h\in
S^2(V^*)}\{\mathcal{R}_h\}$.
\item $\mathcal{R}_{\mathcal{C}}\in\mathfrak{F}(V)$ and
      $\mathfrak{F}(V)=\operatorname{Span}_{\mathcal{C}\in S^2(V^*)\otimes V}\{\mathcal{R}_{\mathcal{C}}\}$.
\end{enumerate}
\end{theorem}

\begin{proof} Assertion (1) follows from the discussion given to establish Theorem \ref{thm-4.1}. We begin the proof of Assertion (2) by
computing:
\begin{eqnarray*}
&&\mathcal{R}_{\mathcal{C}}(v,w)u=\mathcal{C}(w,\mathcal{C}(v,u))-\mathcal{C}(v,\mathcal{C}(w,u))=-\mathcal{R}_{\mathcal{C}}(w,v)u,\\
&&\mathcal{R}_{\mathcal{C}}(w,v)u+\mathcal{R}_{\mathcal{C}}(v,u)w+\mathcal{R}_{\mathcal{C}}(u,w)v
=\mathcal{C}(v,\mathcal{C}(w,u))-\mathcal{C}(w,\mathcal{C}(v,u))\\
&&\qquad+\mathcal{C}(w,\mathcal{C}(u,v))-\mathcal{C}(u,\mathcal{C}(w,v))+\mathcal{C}(u,\mathcal{C}(v,w))-\mathcal{C}(v,\mathcal{C}(u,w))\\
&&\qquad=0\,.
\end{eqnarray*}
Let $C(e_i,e_j)=\sum_kC_{ij}{}^ke_k$ where $\{e_i\}$ is a basis for $V$. We show that $\mathcal{R}_{\mathcal{C}}\in\mathfrak{F}(V)$ by
checking:
\begin{eqnarray*}
&&\mathcal{R}_{\mathcal{C}}(e_i,e_j)e_k=\textstyle\sum_{l,n}\{C_{jl}{}^nC_{ik}{}^l-C_{il}{}^nC_{jk}{}^l\}e_n,\\
&&\rho_{34}(\mathcal{R}_{\mathcal{C}})(e_i,e_j)=\textstyle\sum_{k,l}\{C_{jl}{}^kC_{ik}{}^l-C_{il}{}^kC_{jk}{}^l\}=0\,.
\end{eqnarray*}

Let $\mathfrak{B}(V):=\operatorname{Span}_{\mathcal{C}\in S^2(V^*)\otimes V}\{\mathcal{R}_{\mathcal{C}}\}$. For $\varepsilon\ne0$, let the
non-zero components of $\mathcal{C}$ be given by:
$$
C_{21}{}^1=C_{12}{}^1=C_{11}{}^2=C_{31}{}^1=C_{13}{}^1=C_{11}{}^3=\varepsilon\,.
$$
We have
$$
\rho_{14}(\mathcal{R}_{\mathcal{C}})(e_2,e_2)=\textstyle\sum_{l,i}\{C_{2l}{}^iC_{i2}{}^l-C_{il}{}^iC_{22}{}^l\}
 =\varepsilon^2\ne0\,.
$$
This shows that $\rho_{14}(\mathcal{R}_{\mathcal{C}})(e_2,e_2)\ne0$. Consequently
$$\sigma_{\pi_s\circ\rho_{14}}S^2(V^*)\subset\mathfrak{B}(V)\,.$$
We also compute
\begin{eqnarray*}
&&\mathcal{R}_{\mathcal{C}}(e_1,e_2)e_1=\textstyle\sum_{l,n}\{C_{2l}{}^nC_{11}{}^l-C_{1l}{}^nC_{21}{}^l\}e_n\\
&&\qquad=C_{21}{}^1C_{11}{}^1e_1
        -C_{11}{}^2C_{21}{}^1e_2-C_{11}{}^3C_{21}{}^1e_3\\
&&\qquad=-\varepsilon^2(e_2+e_3)\,.
\end{eqnarray*}
If $\mathcal{R}_{\mathcal{C}}\in\sigma_{\pi_s\circ\rho_{14}}S^2(V^*)$, then
$\mathcal{R}_{\mathcal{C}}(e_1,e_2)e_1\in\operatorname{Span}\{e_1,e_2\}$ which is false. Thus
$$\sigma_{\pi_s\circ\rho_{14}}S^2(V^*)\not\subset\mathfrak{B}(V)\,.$$
The desired result now follows.
\end{proof}

\subsection{The action of $O(V,\langle\cdot,\cdot\rangle)$ on $\mathfrak{F}(V)$}\label{sect-4.2}
We can use Theorems \ref{thm-2.5} and \ref{thm-4.1} to see that there is an $O(V,\langle\cdot,\cdot\rangle)$
equivariant orthogonal decomposition of
\begin{eqnarray*}
\mathfrak{f}(V,\langle\cdot,\cdot\rangle)&\approx&\mathfrak{w}(V,\langle\cdot,\cdot\rangle)\oplus\mathbb{R}\oplus
S_0^2(V^*,\langle\cdot,\cdot\rangle)\\
&\oplus&S_0^2(V^*,\langle\cdot,\cdot\rangle)\oplus\Lambda^2(V^*)\oplus W_7\oplus W_8
\end{eqnarray*}
is a direct sum of $7$ irreducible $O(V,\langle\cdot,\cdot\rangle)$ modules. Since $S_0^2(V^*,\langle\cdot,\cdot\rangle)$
is repeated with multiplicity 2, the decomposition is not unique.

\par We now make this decomposition a bit more explicit  to identify the factors $W_7$ and $W_8$. We adopt the notation of
Equation (\ref{eqn-3.a}) and let $\id\otimes\pi_s$ symmetrize the last two components of
$T\in\otimes^4V^*$. Let $\sigma_{\id\otimes\pi_s}$ be the splitting of Equation (\ref{eqn-3.c}). Finally,
let $\alpha$ be the map of Equation (\ref{eqn-3.d}).

\begin{lemma}\label{lem-4.4}
We  have an $O(V,\langle\cdot,\cdot\rangle)$ equivariant short exact sequence
$$0\rightarrow\mathfrak{a}(V)\rightarrow\mathfrak{f}(V,\langle\cdot,\cdot\rangle)\mapright{\id\otimes\pi_s}
\Lambda^2(V^*)\otimes S^2_0(V^*,\langle\cdot,\cdot\rangle)\rightarrow0$$
which is equivariantly split by the map $\sigma_{\id\otimes\pi_s}$.
\end{lemma}

\begin{proof} Let $F\in\mathfrak{f}(V,\langle\cdot,\cdot\rangle)$. We have 
\begin{eqnarray*}
&&(\id\otimes\pi_s)(F)(x,y,z,w)=\textstyle\frac12\{F(x,y,z,w)+F(x,y,w,z)\},\\
&&(\id\otimes\pi_s)(F)=0\quad\Leftrightarrow\quad
F(x,y,z,w)=-F(x,y,w,z)\ \forall\ x,y,z,w\in V\,.
\end{eqnarray*}
This implies $F\in\mathfrak{a}(V)$. Conversely, if $F\in\mathfrak{a}(V)$,
then
$\rho_{34}(F)=0$ and $(\id\otimes\pi_s)F=0$ and hence $F\in\mathfrak{f}(V,\langle\cdot,\cdot\rangle)$. Thus
$$\ker\{\id\otimes\pi_s\}\cap\mathfrak{f}(V,\langle\cdot,\cdot\rangle)=\mathfrak{a}(V)\,.$$
Furthermore 
$$\rho_{34}(F)=(\id\otimes\Tr)((\id\otimes\pi_s)F)$$
and consequently $(\id\otimes\pi_s)$ takes values in $\Lambda^2(V^*)\otimes S_0^2(V^*,\langle\cdot,\cdot\rangle)$.

In the proof of Theorem \ref{thm-3.2}, we showed that $\sigma_{\id\otimes\pi_s}$ takes values in $\mathfrak{r}(V)$ and that
$(\id\otimes\pi_s)\sigma_{\id\otimes\pi_s}$ is the identity on $\Lambda^2(V)\otimes S^2(V^*)$. Thus
$\sigma_{\id\otimes\pi_s}S\in\mathfrak{f}(V,\langle\cdot,\cdot\rangle)$ if and only if $S\in\Lambda^2(V^*)\otimes
S^2_0(V^*,\langle\cdot,\cdot\rangle)$. \end{proof}

This shows that
\begin{eqnarray*}
&&\mathfrak{f}(V,\langle\cdot,\cdot\rangle)\approx\mathfrak{a}(V)\oplus\Lambda^2(V^*)\otimes
S_0^2(V^*,\langle\cdot,\cdot\rangle),\quad\text{so}\\ 
&&\Lambda^2(V^*)\otimes S_0^2(V^*,\langle\cdot,\cdot\rangle)\approx
S_0^2(V^*,\langle\cdot,\cdot\rangle)\oplus\Lambda^2(V^*)\oplus W_7\oplus W_8\,.
\end{eqnarray*}
We therefore study
$\Lambda^2(V^*)\otimes S^2(V^*)$ as an
$O(V,\langle\cdot,\cdot\rangle)$ module and identify the copies of $\Lambda^2(V^*)$ and $S_0^2(V^*,\langle\cdot,\cdot\rangle)$ in
$\Lambda^2(V^*)\otimes S_0^2(V^*,\langle\cdot,\cdot\rangle)$. Let 
$$\Theta\in\Lambda^2(V^*)\otimes S_0^2(V^*,\langle\cdot,\cdot\rangle),\quad
\psi\in S_0^2(V^*,\langle\cdot,\cdot\rangle),\quad
\omega\in\Lambda^2(V^*)\,.
$$
Let $\Xi$ be as in Equation (\ref{eqn-1.r}). Define:\def\gronka{\textstyle\frac
m{m^2-4}}\def\gronkb{\textstyle\frac4m}
\medbreak\qquad\qquad
$\pi_{1,s}(\Theta)_{jk}:=(\pi_s(\rho_{14}\theta))_{jk}=\textstyle\frac12\sum_{il}\Xi^{il}\{\Theta_{ijkl}+\Theta_{ikjl}\}$,
\smallbreak\qquad\qquad
$\pi_{1,a}(\Theta)_{jk}:=(\pi_a(\rho_{14}\theta))_{jk}=\textstyle\frac12\sum_{il}\Xi^{il}\{\Theta_{ijkl}-\Theta_{ikjl}\}$,
\smallbreak\qquad\qquad
$\pi_\Lambda(\Theta)_{ijkl}:=\textstyle\frac12(\Theta_{kjil}+\Theta_{ikjl}-\Theta_{ljik}-\Theta_{iljk})$,
\smallbreak\qquad\qquad
$\sigma_{\pi_{1,s}}(\psi)_{ijkl}:=\textstyle\frac1m\{\Xi_{il}\psi_{jk}-\Xi_{jl}\psi_{ik}+\Xi_{ik}\psi_{jl}-\Xi_{jk}\psi_{il}\}$,
\smallbreak\qquad\qquad
$\sigma_{\pi_{1,a}}(\omega)_{ijkl}:=\gronka\{\Xi_{il}\omega_{jk}+\Xi_{ik}\omega_{jl}-\Xi_{jl}\omega_{ik}-\Xi_{jk}\omega_{il}
+\gronkb\omega_{ij}\Xi_{kl}\}$,
\smallbreak\qquad\qquad
$\sigma_{\pi_\Lambda}(\Theta)_{ijkl}:=\textstyle\frac12(\Theta_{kjil}-\Theta_{kijl})$,
\smallbreak\qquad\qquad
$\Lambda_0^2(\Lambda^2(V^*)):=\left\{\Theta:\Theta_{ijkl}=-\Theta_{jikl}=-\Theta_{klij},\quad
  \textstyle\sum_{il}\Xi^{il}\Theta_{ijkl}=0\right\}$,
\smallbreak\qquad\qquad
$\mathfrak{S}(V,\langle\cdot,\cdot\rangle):=\ker\{\pi_{1,s}\}\cap\ker\{\pi_{1,a}\}
   \cap\ker\{\pi_\Lambda\}\cap\Lambda^2(V^*)
\otimes S^2_0(V^*,\langle\cdot,\cdot\rangle)$.

\begin{lemma}\label{lem-4.5}
We have $O(V,\langle\cdot,\cdot\rangle)$ equivariant short exact sequences
\begin{eqnarray*}
&&0\rightarrow\ker\{\pi_{1,s}\}\rightarrow\Lambda^2(V^*)\otimes S_0^2(V^*,\langle\cdot,\cdot\rangle)\mapright{\pi_{1,s}}
 S_0^2(V^*,\langle\cdot,\cdot\rangle)\rightarrow0,\\
&&0\rightarrow\ker\{\pi_{1,a}\}\rightarrow\Lambda^2(V^*)\otimes S_0^2(V^*,\langle\cdot,\cdot\rangle)
\mapright{\pi_{1,a}}\Lambda^2(V^*)\rightarrow0,\\
&&0\rightarrow\ker\{\pi_{1,a}\}\cap\ker\{\pi_\Lambda\}
\rightarrow\ker\{\pi_{1,a}\}\mapright{\pi_\Lambda}\Lambda_0^2(\Lambda^2(V^*))\rightarrow0
\,.
\end{eqnarray*}
These sequences are equivariantly split, respectively, by $\sigma_{\pi_{1,s}}$, $\sigma_{\pi_{1,a}}$,
and
$\sigma_{\pi_\Lambda}$. This gives an $O(V,\langle\cdot,\cdot\rangle)$ equivariant decomposition of
$$\Lambda^2(V^*)\otimes S_0^2(V^*,\langle\cdot,\cdot\rangle)\approx S_0^2(V^*,\langle\cdot,\cdot\rangle)
\oplus\Lambda^2(V^*)\oplus\Lambda_0^2(\Lambda^2(V^*))\oplus\mathfrak{S}(V,\langle\cdot,\cdot\rangle)$$
as the direct sum of irreducible $O(V,\langle\cdot,\cdot\rangle)$ modules where
$$\begin{array}{l}
\dim\{S_0^2(V^*,\langle\cdot,\cdot\rangle)\}=\textstyle\frac{(m-1)(m+2)}2,\\
\dim\{\Lambda^2(V^*)\}=\textstyle\frac{m(m-1)}2,\vphantom{\vrule height 11pt}\\
\dim\{\Lambda_0^2(\Lambda^2(V^*))\}
   =\frac{m(m-1)(m-3)(m+2)}8,\vphantom{\vrule height 11pt}\\
\dim\{\mathfrak{S}(V,\langle\cdot,\cdot\rangle)\}
   =\textstyle\frac{(m-1)(m-2)(m+1)(m+4)}8,\vphantom{\vrule height 11pt}\\
\dim\{\Lambda^2(V^*)\otimes S_0^2(V^*,\langle\cdot,\cdot\rangle)\}=\textstyle\frac{m(m-1)^2(m+2)}8\,.
   \vphantom{\vrule height 11pt}
\end{array}$$
We have $W_8\approx\Lambda_0^2(\Lambda^2(V^*))$ and $W_7\approx\mathfrak{S}(V,\langle\cdot,\cdot\rangle)$.
\end{lemma}

\begin{proof} It is clear that $\pi_{1,s}$ takes values in $S^2(V^*)$. Let $\Xi$ be as in Equation (\ref{eqn-1.r}). We
show that $\pi_{1,s}$ takes values in $S^2_0(V^*,\langle\cdot,\cdot\rangle)$ by checking:
\begin{eqnarray*}
\Tr\{\pi_{1,s}(\Theta)\}&=&\textstyle\frac12\sum_{ijkl}\Xi^{il}\Xi^{jk}\{\Theta_{ijkl}+\Theta_{ikjl}\}\\
&=&\textstyle\sum_{ijkl}\Xi^{il}\Xi^{jk}\Theta_{ijkl}=\sum_{ijkl}\Xi^{jk}\Xi^{il}\Theta_{jilk}\\
&=&\textstyle
-\sum_{ijkl}\Xi^{jk}\Xi^{il}\Theta_{ijkl}=-\Tr\{\pi_{1,s}(\Theta)\}\,.
\end{eqnarray*}
It is clear that $\sigma_{\pi_{1,s}}$ takes values in $\Lambda^2(V^*)\otimes S^2(V^*)$. We verify that
$\sigma_{\pi_{1,s}}$ takes values in $\Lambda^2(V^*)\otimes S^2_0(V^*,\langle\cdot,\cdot\rangle)$ by checking the trace
condition:
\begin{eqnarray*}
\textstyle\sum_{kl}\Xi^{kl}\sigma_{\pi_{1,s}}(\psi)_{ijkl}
&=&\textstyle\frac1m\textstyle\sum_{kl}\Xi^{kl}\{\Xi_{il}\psi_{jk}-\Xi_{jl}\psi_{ik}+\Xi_{ik}\psi_{jl}-\Xi_{jk}\psi_{il}\}\\
&=&\textstyle\frac1m\{\psi_{ji}-\psi_{ij}+\psi_{ji}-\psi_{ij}\}=0\,.
\end{eqnarray*}
We check that $\sigma_{\pi_{1,s}}$ is a splitting by verifying:
\begin{eqnarray*}
\pi_{1,s}(\sigma_{\pi_{1,s}}(\psi))_{jk}
&=&\textstyle\frac1{2m}\sum_{il}\Xi^{il}\{\Xi_{il}\psi_{jk}-\Xi_{jl}\psi_{ik}+\Xi_{ik}\psi_{jl}-\Xi_{jk}\psi_{il}\\
&&\qquad\quad\phantom{.}+\Xi_{il}\psi_{kj}-\Xi_{kl}\psi_{ij}+\Xi_{ij}\psi_{kl}-\Xi_{kj}\psi_{il}\}\\
&=&\textstyle\frac1{2m}\{m\psi_{jk}-\psi_{jk}+\psi_{jk}-\Xi_{jk}\Tr\{\psi\}\\
&&\quad+m\psi_{kj}-\psi_{kj}+\psi_{kj}-\Xi_{kj}\Tr\{\psi\}\}\\
&=&\psi_{jk}\,.
\end{eqnarray*}

Clearly $\pi_{1,a}$ takes values in $\Lambda^2(V^*)$ and $\sigma_{\pi_{1,a}}$ takes values in
$\Lambda^2(V^*)\otimes S^2(V^*)$. We check the trace condition by computing:
\begin{eqnarray*}
&&\{(\id\otimes\Tr)(\sigma_{\pi_{1,a}}(\omega))\}_{ij}\\
&=&\gronka\sum_{kl}\Xi^{kl}\{\Xi_{il}\omega_{jk}+\Xi_{ik}\omega_{jl}
-\Xi_{jl}\omega_{ik}-\Xi_{jk}\omega_{il}+\gronkb\omega_{ij}\Xi_{kl}\}\\
&=&\gronka\{\omega_{ji}+\omega_{ji}-\omega_{ij}-\omega_{ij}+\gronkb m\omega_{ij}\}\\
&=&\gronka(-4+\gronkb m)\omega_{ij}
   =0\,.
\end{eqnarray*}
To check $\sigma_{\pi_{1,a}}$ is a splitting, we compute:
\begin{eqnarray*}
&&\pi_{1,a}(\sigma_{\pi_{1,a}}(\omega))_{jk}\\
&=&\textstyle\frac12\gronka\sum_{il}\Xi^{il}\{\Xi_{il}\omega_{jk}+\Xi_{ik}\omega_{jl}-\Xi_{jl}\omega_{ik}-\Xi_{jk}\omega_{il}
   +\gronkb\omega_{ij}\Xi_{kl}\\
&-&\Xi_{il}\omega_{kj}-\Xi_{ij}\omega_{kl}+\Xi_{kl}\omega_{ij}+\Xi_{kj}\omega_{il}-\gronkb\omega_{ik}\Xi_{jl}\}\\
&=&\textstyle\frac12\gronka\{m\omega_{jk}+\omega_{jk}-\omega_{jk}+\gronkb\omega_{kj}
   -m\omega_{kj}-\omega_{kj}+\omega_{kj}-\gronkb\omega_{jk}\}\\
&=&\gronka\{m-\gronkb\}\omega_{jk}=\omega_{jk}\,.
\end{eqnarray*}

Let $S\in\ker\{\pi_{1,a}\}\cap\{\Lambda^2(V^*)\otimes S_0^2(V^*,\langle\cdot,\cdot\rangle)\}$. To check that
$\pi_\Lambda$ takes values in
$\Lambda_0^2(\Lambda^2(V^*))$, we compute:
\begin{eqnarray*}
&&\pi_\Lambda(S)_{ijkl}=\textstyle\frac12(S_{kjil}+S_{ikjl}-S_{ljik}-S_{iljk}),\\
&&\pi_\Lambda(S)_{jikl}=\textstyle\frac12(S_{kijl}+S_{jkil}-S_{lijk}-S_{jlik})=-\pi_\Lambda(S)_{ijkl},\\
&&\pi_\Lambda(S)_{klij}=\textstyle\frac12(S_{ilkj}+S_{kilj}-S_{jlki}-S_{kjli})=-\pi_\Lambda(S)_{ijkl},\\
&&\rho_{14}(\pi_\Lambda(S))_{jk}=\textstyle\frac12\sum_{il}\Xi^{il}\{S_{kjil}+S_{ikjl}-S_{ljik}-S_{iljk}\}\\
&&\qquad\qquad=\{\textstyle\frac12\rho_{34}(S)+\pi_{1,a}(S)\}_{jk}=0\,.
\end{eqnarray*}
Let  $T\in\Lambda_0^2(\Lambda^2(V^*))$. To check $\sigma_{\pi_\Lambda}$ takes values in $\Lambda^2(V^*)\otimes
S_0^2(V^*,\langle\cdot,\cdot\rangle)$, we compute:
\begin{eqnarray*}
&&\sigma_{\pi_\Lambda}(T)_{ijkl}=\textstyle\frac12(T_{kjil}-T_{kijl}),\\
&&\sigma_{\pi_\Lambda}(T)_{jikl}=\textstyle\frac12(T_{kijl}-T_{kjil})=-\sigma_{\pi_\Lambda}(T)_{ijkl},\\
&&\sigma_{\pi_\Lambda}(T)_{ijlk}=\textstyle\frac12(T_{ljik}-T_{lijk})=\frac12(T_{jkli}-T_{iklj})\\
&&\qquad=\textstyle\frac12(T_{kjil}-T_{kijl})=\sigma_{\pi_\Lambda}(T)_{ijlk},\\
&&\textstyle\sum_{kl}\Xi^{kl}\sigma_{\pi_\Lambda}(T)_{ijkl}=\frac12\sum_{kl}\Xi^{kl}(T_{kjil}-T_{kijl})=0\,.
\end{eqnarray*}
Finally, we verify that $\sigma_{\pi_\Lambda}$ is a splitting by computing
\begin{eqnarray*}
&&\{\pi_\Lambda(\sigma_{\pi_\Lambda}(T))\}_{ijkl}\\
&=&\textstyle\frac12(\sigma_{\pi_\Lambda}(T)_{kjil}+\sigma_{\pi_\Lambda}(T)_{ikjl}
     -\sigma_{\pi_\Lambda}(T)_{ljki}-\sigma_{\pi_\Lambda}(T)_{ilkj})\\
&=&\textstyle\frac14(T_{ijkl}-T_{ikjl}+T_{jkil}-T_{jikl}
-T_{kjli}+T_{klji}-T_{klij}+T_{kilj})\\
&=&T_{ijkl}\,.
\end{eqnarray*}

We compute dimensions:
\begin{eqnarray*}
&&\dim\{\Lambda^2(V^*)\}=\textstyle\frac12m(m-1),\\
&&\dim\{\Lambda^2(\Lambda^2(V^*))\}=\textstyle\frac12\{\frac12m(m-1)\}\{\frac12m(m-1)-1\},\\
&&\dim\{\Lambda_0^2(\Lambda^2(V^*))\}=\dim\{\Lambda^2(\Lambda^2(V^*))\}-\dim\{\Lambda^2(V^*)\}\\
&=&\textstyle\frac12\{\frac12m(m-1)\}\{\frac12m(m-1)-1\}-\frac12m(m-1)\\
&=&\textstyle\{\frac12m(m-1)\}\{\frac14m(m-1)-\frac12-1\}\\
&=&\textstyle\frac18\{m(m-1)\}\{m(m-1)-6\}=\frac18m(m-1)(m-3)(m+2)\\
&=&\dim\{W_8\}
\end{eqnarray*}
and
\begin{eqnarray*}
&&\dim\{\mathfrak{S}(V,\langle\cdot,\cdot\rangle)\}\\
&=&\dim\{\Lambda^2(V^*)\otimes S_0^2(V^*,\langle\cdot,\cdot\rangle)\}-\dim\{\Lambda_0^2(\Lambda^2(V^*))\}\\
&-&\dim\{S_0^2(V^*,\langle\cdot,\cdot\rangle)\}-\dim\{\Lambda^2(V^*)\}\\
&=&\dim\{\Lambda^2(V^*)\otimes
S_0^2(V^*,\langle\cdot,\cdot\rangle)\}-\dim\{\Lambda^2(\Lambda^2(V^*))\}-\dim\{S_0^2(V^*,\langle\cdot,\cdot\rangle\}\\
&=&\textstyle\frac{m(m-1)(m-1)(m+2)}4-\textstyle\frac{m(m-1)(m(m-1)-2)}8-\textstyle\frac{(m-1)(m+2)}2\\
&=&\textstyle\frac{m-1}8\{2m(m-1)(m+2)-m(m-2)(m+1)-4(m+2)\}\\ &=&\textstyle\frac{(m-1)(m-2)(m+1)(m+4)}8=\dim\{W_7\}\,.
\end{eqnarray*}
The remaining assertions now follow from Theorem \ref{thm-2.5} (1); this also establishes Theorem \ref{thm-2.5} (2c).
\end{proof}

As an immediate consequence, we have
\begin{theorem}\label{thm-4.6}
\ \begin{enumerate}
\item There is an
$O(V,\langle\cdot,\cdot\rangle)$ equivariant orthogonal decomposition of
$$\mathfrak{F}(V)\approx\mathfrak{f}(V)=W_1\oplus W_2\oplus W_4\oplus W_5\oplus W_6\oplus W_7\oplus W_8$$
as the direct sum of irreducible $O(V,\langle\cdot,\cdot\rangle)$ modules where:
$$\begin{array}{ll}
\dim\{W_1\}=1,&
\dim\{W_2\}=\dim\{W_5\}=\textstyle\frac{(m-1)(m+2)}2,\\
\dim\{W_4\}=\textstyle\frac{m(m-1)}2,&
\dim\{W_6\}=\textstyle\frac{m(m+1)(m-3)(m+2)}{12},\\
\dim\{W_7\}=\textstyle\frac{(m-1)(m-2)(m+1)(m+4)}8,&
\dim\{W_8\}=\textstyle\frac{m(m-1)(m-3)(m+2)}8\,.
\end{array}$$
\item There are the following isomorphisms as $O(\langle\cdot,\cdot\rangle)$
modules:
\begin{enumerate}
\item $W_1\approx\mathbb{R}$,
  $W_2\approx W_5\approx S_0^2(V^*,\langle\cdot,\cdot\rangle)$, and
  $W_4\approx\Lambda^2(V^*)$.
\item $W_6\approx\mathfrak{w}(V,\langle\cdot,\cdot\rangle)$ is the space of Weyl conformal curvature tensors.
\item  $W_7\approx\mathfrak{S}(V,\langle\cdot,\cdot\rangle)$ and $W_8\approx\Lambda_0^2(\Lambda^2(V^*))$.
\end{enumerate}
\end{enumerate}
\end{theorem}

\section{The proof of Theorem \ref{thm-3.1}}\label{sect-5}

Let $\mathfrak{b}$ be a non-empty subspace of $\mathfrak{a}(V)$ which is invariant under the action of $GL(V)$. We must show
that $\mathfrak{b}=\mathfrak{a}(V)$. Choose a positive definite inner product
$\langle\cdot,\cdot\rangle$ on $V$. Then
$\mathfrak{b}$ is invariant under the action of $O(V,\langle\cdot,\cdot\rangle)$ as well.  Let
$\pi_{\mathfrak{w}}$, $\pi_0$, and $\pi_{\mathbb{R}}$ be the projections on the appropriate module
summands in the decomposition of Theorem \ref{thm-3.2} (3);
\begin{eqnarray*}
&&\pi_{\mathbb{R}}(R):=\tau(\rho_{14}(R)),\quad
\pi_0(R):=\rho_{14}(R)-\textstyle\frac1m\tau(\rho_{14}(R))\langle\cdot,\cdot\rangle,\\
&&\pi_{\mathfrak{w}}(R):=R-\sigma_{\mathfrak{a},\rho_{14}}(\rho_{14}(R))\quad\text{where}\\
&&\sigma_{\mathfrak{a},\rho_{14}}(\psi):=
   \textstyle\frac2{m-2}\psi\wedge\langle\cdot,\cdot\rangle-\frac{\tau(\psi)}{(m-1)(m-2)}
 \langle\cdot,\cdot\rangle\wedge\langle\cdot,\cdot\rangle\,.
\end{eqnarray*}
 Since $O(V,\langle\cdot,\cdot\rangle)$ is
a compact Lie group acting orthogonally, the projections are orthogonal projections. Furthermore:
$$\begin{array}{lll}
\pi_{\mathfrak{w}}(\mathfrak{b})\ne\{0\}&\Rightarrow&\mathfrak{w}(V,\langle\cdot,\cdot\rangle)\subset\mathfrak{b},\\
\pi_0(\mathfrak{b})\ne\{0\}&\Rightarrow&\sigma_{\mathfrak{a},\rho_{14}}(S_0^2(V,\langle\cdot,\cdot\rangle))\subset\mathfrak{b},\\
\pi_{\mathbb{R}}(\mathfrak{b})\ne\{0\}&\Rightarrow&\sigma_{\mathfrak{a},\rho_{14}}(\langle\cdot,\cdot\rangle)
   \subset\mathfrak{b}\,.
\end{array}$$
Let $\{e_i\}$ be an orthonormal basis for $V$. Let $\{\lambda_i\}$ be distinct positive constants. Define
$\Theta\in GL(V)$ by setting:
$$\Theta(e_i)=\lambda_ie_i\,.$$

Suppose $\pi_{\mathbb{R}}(\mathfrak{b})\ne\{0\}$. The component corresponding to $\mathbb{R}$ in $\mathfrak{a}(V)$
is generated by $A:=\langle\cdot,\cdot\rangle\wedge\langle\cdot,\cdot\rangle$. Consequently $A\in\mathfrak{b}$; the non-zero
components of
$\Theta^*(A)$ and $\rho_{14}(\Theta^*(A))$ are, up to the usual $\mathbb{Z}_2$ symmetries and modulo a suitable normalizing
constant which plays no role, given by
\begin{eqnarray*}
&&\Theta^*(A)(e_i,e_j,e_j,e_i)=\lambda_i^2\lambda_j^2\quad\text{and}\quad
  \rho_{14}(\Theta^*(A))(e_i,e_i)=\lambda_i^2\textstyle\sum_{j\ne i}\lambda_j^2\,.
\end{eqnarray*}
This shows the projection of $\Theta^*(A)$, and hence of $\mathfrak{b}$, on $S_0(V^*,\langle\cdot,\cdot\rangle)$ is non-zero.
Let $A_1$ be the algebraic curvature tensor whose only non-zero component, up to the usual $\mathbb{Z}_2$ symmetries, is
$A_1(e_1,e_2,e_2,e_1)$. As
$\mathfrak{b}$ is closed, we show that $A_1\in\mathfrak{b}$ by taking the limit
$$\lambda_1\rightarrow1,\quad\lambda_2\rightarrow1,\quad\lambda_j\rightarrow0\text{ for }j\ge3\,.$$ 
As
$\{A_1-\sigma_{\mathfrak{a},\rho_{14}}(\rho_{14}(A_1))\}(e_1,e_3,e_3,e_1)\ne0$, one has
$\pi_{\mathfrak{w}}(\mathfrak{b})\ne0$. We summarize:
$$\pi_{\mathbb{R}}(\mathfrak{b})\ne0\quad\Rightarrow\quad\mathfrak{b}=\mathfrak{a}(V)\,.$$

Suppose $\pi_0(\mathfrak{b})\ne0$. Then
$\sigma_{\mathfrak{a},\rho_{14}}(S_0^2(V^*,\langle\cdot,\cdot\rangle))\subset\mathfrak{b}$. Define
$\psi\in S_0^2(V^*,\langle\cdot,\cdot\rangle)$ with  non-zero components
$$\psi(e_1,e_1)=1,\quad\psi(e_2,e_2)=1,\quad\text{and}\quad\psi(e_3,e_3)=-2\,.$$
Let $A=\sigma_{\mathfrak{a},\rho_{14}}(\psi)=\textstyle\frac2{m-2}\psi\wedge\langle\cdot,\cdot\rangle\in\mathfrak{b}$. We compute:
\begin{eqnarray*}
&&\Theta^*(A)(e_i,e_j,e_k,e_l)=\lambda_i\lambda_j\lambda_k\lambda_l\textstyle\frac2{m-2}
   \{\psi(e_i,e_l)\delta_{jk}+\psi(e_j,e_k)\delta_{il}\\
&&\qquad\qquad-\psi(e_i,e_k)\delta_{jl}-\psi(e_j,e_l)\delta_{ik}\},\\
&&\tau(\rho_{14}(\Theta^*(A)))=\textstyle\frac2{m-2}\sum_{i,j}\Theta^*A(e_i,e_j,e_j,e_i)\\
&=&\textstyle\frac2{m-2}\{(\lambda_1^2+\lambda_2^2-2\lambda_3^2)(\sum_j\lambda_j^2)\}-
  \textstyle\frac2{m-2}\sum_i\{\lambda_1^4+\lambda_2^4-2\lambda_3^4\}\,.
\end{eqnarray*}
This is non-zero for generic values of $\vec\lambda$. This shows
$\pi_{\mathbb{R}}(\mathfrak{b})\ne\{0\}$. Combining this result with the result of the
previous paragraph yields:
$$\rho_{14}(\mathfrak{b})\ne\{0\}\quad\Rightarrow\quad\mathfrak{b}=\mathfrak{a}(V)\,.$$

Finally, suppose $\pi_{\mathfrak{w}}(\mathfrak{b})\ne0$. Then $\mathfrak{w}(V,\langle\cdot,\cdot\rangle)\subset\mathfrak{b}$.
Let $A\in\mathfrak{a}$ be defined with non-zero components, up to the usual $\mathbb{Z}_2$
symmetries, by 
$$A(e_1,e_3,e_4,e_1)=+1\quad\text{and}\quad A(e_2,e_3,e_4,e_2)=-1\,.$$
Then $\rho_{14}(A)=0$ so $A\in\mathfrak{w}(V,\langle\cdot,\cdot\rangle)$. We have
$$\Theta^*(A)(e_1,e_3,e_4,e_1)=\lambda_1^2\lambda_3\lambda_4\quad\text{and}\quad
  \Theta^*(A)(e_2,e_3,e_4,e_2)=\lambda_2^2\lambda_3\lambda_4\,.$$
Thus $\rho_{14}(\Theta^*(A))(e_3,e_4)=\lambda_3\lambda_4(\lambda_1^2-\lambda_2^2)\ne0$. Since
$\rho_{14}(\Theta^*(A))\ne0$ we may conclude that
$\mathfrak{b}=\mathfrak{a}(V)$.
\hfill\qedbox

\section*{Acknowledgments} Research of N.Bla\v zi\'c partially supported by the DAAD
 (Germany), TU Berlin, and MNTS 1854 (Serbia).  Research of 
 P. Gilkey partially supported by the Max Planck Institute in the 
 Mathematical Sciences (Leipzig, Germany) and by DFG PI158/4-4.
 Research of S. Nik\v cevi\'c partially supported by DAAD
 (Germany), TU Berlin, MM 1646 (Serbia), and the Dierks von Zweck-Stiftung.  Research of 
 U. Simon partially supported by DFG PI158/4-4.

\end{document}